\documentclass[12pt,a4paper,draft]{amsart}

\usepackage{latexsym}
\usepackage{ifthen}
\usepackage[leqno]{amsmath}
\usepackage{enumerate}
\usepackage{calc}
\usepackage{amstext,amsbsy,amsopn,amsthm,amsgen,amsfonts,amscd,amsxtra,upref}
\usepackage{mathrsfs}\usepackage{euscript}\usepackage{amssymb}

\swapnumbers
\theoremstyle{plain}
\newtheorem{Thm}{Theorem}[section]
\newtheorem{Lem}[Thm]{Lemma}
\newtheorem{Cor}[Thm]{Corollary}
\newtheorem{Pro}[Thm]{Proposition}
\newtheorem{Prp}[Thm]{Properties}
\newtheorem{Sub}[Thm]{Sublemma}

\theoremstyle{definition}
\newtheorem{Def}[Thm]{Definition}
\newtheorem{Exm}[Thm]{Example}
\newtheorem{Exs}[Thm]{Examples}

\theoremstyle{remark}
\newtheorem{Rem}[Thm]{Remark}
\newtheorem{Rms}[Thm]{Remarks}
\newtheorem*{Com}{Commentary}

\newcommand{\myEmail}{piotr.niemiec@uj.edu.pl}
\newcommand{\myAddress}{\noindent{}Piotr Niemiec\\{}Jagiellonian University\\{}Institute of Mathematics\\{}
   ul. \L{}ojasiewicza 6\\{}30-348 Krak\'{o}w\\{}Poland}
\newcommand{\myData}{\author[P. Niemiec]{Piotr Niemiec}\address{\myAddress}\email{\myEmail}}

\newcommand{\NNN}{\mathbb{N}}
\newcommand{\RRR}{\mathbb{R}}

\newcommand{\CCc}{\CMcal{C}}
\newcommand{\FFf}{\CMcal{F}}\newcommand{\GGg}{\CMcal{G}}

\newcommand{\SSs}{\CMcal{S}}
\newcommand{\UUu}{\CMcal{U}}\newcommand{\VVv}{\CMcal{V}}\newcommand{\WWw}{\CMcal{W}}

\newcommand{\SsS}{\EuScript{S}}

\newcommand{\mM}{\mathfrak{m}}

\newcommand{\SECT}[1]{\section{#1}\renewcommand{\theequation}{\thesection-\arabic{equation}}\setcounter{equation}{0}}

\newcounter{help}
\newcommand{\ITE}[3]{\ifthenelse{#1}{#2}{#3}}\newcommand{\ITEE}[3]{\ITE{\equal{#1}{#2}}{#3}{}}

\newcommand{\dist}{\operatorname{dist}}

\newcommand{\im}{\operatorname{im}}\newcommand{\id}{\operatorname{id}}
\newcommand{\cll}{\operatorname{cl}}\newcommand{\intt}{\operatorname{int}}

\newcommand{\Metr}{\operatorname{Metr}}
\newcommand{\cov}{\operatorname{cov}}
\newcommand{\Emb}{\operatorname{Emb}}

\newcommand{\leqsl}{\leqslant}\newcommand{\geqsl}{\geqslant}

\newcommand{\epsi}{\varepsilon}\newcommand{\varempty}{\varnothing}\newcommand{\dd}{\colon}

\newcommand{\THM}[1]{Theorem~\textup{\ref{thm:#1}}}
\newcommand{\COR}[1]{Corollary~\textup{\ref{cor:#1}}}
\newcommand{\LEM}[1]{Lemma~\textup{\ref{lem:#1}}}\newcommand{\PRO}[1]{Proposition~\textup{\ref{pro:#1}}}

\newcommand{\EXM}[1]{Example~\textup{\ref{exm:#1}}}

\newenvironment{thm}[1]{\begin{Thm}\label{thm:#1}}{\end{Thm}}\newenvironment{lem}[1]{\begin{Lem}\label{lem:#1}}{\end{Lem}}
\newenvironment{cor}[1]{\begin{Cor}\label{cor:#1}}{\end{Cor}}\newenvironment{pro}[1]{\begin{Pro}\label{pro:#1}}{\end{Pro}}
\newenvironment{dfn}[1]{\begin{Def}\label{def:#1}}{\end{Def}}
\newenvironment{exm}[1]{\begin{Exm}\label{exm:#1}}{\end{Exm}}
\newenvironment{rem}[1]{\begin{Rem}\label{rem:#1}}{\end{Rem}}

\newenvironment{thmm}[2]{\begin{Thm}[#2]\label{thm:#1}}{\end{Thm}}
\newenvironment{lemm}[2]{\begin{Lem}[#2]\label{lem:#1}}{\end{Lem}}
\newenvironment{corr}[2]{\begin{Cor}[#2]\label{cor:#1}}{\end{Cor}}

\newcommand{\bibITEM}[2]{\ITE{\equal{#2}{}}{\bibitem{#1} }{\bibitem[#2]{#1} }}
\newcommand{\BIB}[8]{
   \bibITEM{#1}{#8} #2, \textit{#3}, #4{} \textbf{#5} (#6), #7.}
\newcommand{\myBIB}[6][P. Niemiec]{#1, \textit{#2}, #3{}\ITE{\equal{#4}{}}{}{ \textbf{#4}} (#5), #6.}
\newcommand{\BIb}[6]{
   \bibITEM{#1}{#6} #2, \textit{#3}, #4, #5.}
\newcommand{\BiB}[9]{
   \bibITEM{#1}{#9} #2, \textit{#3}, #4{} \textit{#5}, #6, #7, #8.}

\newcommand{\myBAPP}[3][P. Niemiec]{
   #1, \textit{#2}, #3}

\newcommand{\jRN}[2][]{
   \ITEE{#2}{ActaM}{\ITE{\equal{#1}{+}}
      {Acta Mathematica}{Acta Math.}}
   \ITEE{#2}{ActaMSinES}{\ITE{\equal{#1}{+}}
      {Acta Mathematica Sinica (English Series)}{Acta Math. Sin. (Engl. Ser.)}}
   \ITEE{#2}{AdvM}{\ITE{\equal{#1}{+}}
      {Advances in Mathematics}{Adv. in Math.}}
   \ITEE{#2}{ACS}{\ITE{\equal{#1}{+}}
      {Applied Categorical Structures}{Appl. Categor. Struct.}}
   \ITEE{#2}{ActaSM}{\ITE{\equal{#1}{+}}
      {Acta Scientiarum Mathematicarum}{Acta Sci. Math.}}
   \ITEE{#2}{AmJM}{\ITE{\equal{#1}{+}}
      {American Journal of Mathematics}{Amer. J. Math.}}
   \ITEE{#2}{AmMMon}{\ITE{\equal{#1}{+}}
      {American Mathematical Monthly}{Amer. Math. Mon.}}
   \ITEE{#2}{AnnSciEcNormSupT}{\ITE{\equal{#1}{+}}
      {Annales Scientifiques de l'\'{E}cole Normale Sup\'{e}rieure (3)}{Ann. Sci. \'{E}c. Norm. Sup\'{e}r. (3)}}
   \ITEE{#2}{AnnM}{\ITE{\equal{#1}{+}}
      {Annals of Mathematics}{Ann. Math.}}
   \ITEE{#2}{AnnProb}{\ITE{\equal{#1}{+}}
      {The Annals of Probability}{Ann. Probab.}}
   \ITEE{#2}{AnnPALog}{\ITE{\equal{#1}{+}}
      {Annals of Pure and Applied Logic}{Ann. Pure Appl. Logic}}
   \ITEE{#2}{APM}{\ITE{\equal{#1}{+}}
      {Annales Polonici Mathematici}{Ann. Polon. Math.}}
   \ITEE{#2}{ArchM}{\ITE{\equal{#1}{+}}
      {Archiv der Mathematik}{Arch. Math.}}
   \ITEE{#2}{AttiAccLincRendNat}{\ITE{\equal{#1}{+}}
      {Atti della Accademia Nazionale dei Lincei. Rendiconti. Classe di Scienze Fisiche, Matematiche e Naturali}
      {Atti Accad. Naz. Lincei Rend. Cl. Sci. Fis. Mat. Nat.}}
   \ITEE{#2}{BAMS}{\ITE{\equal{#1}{+}}
      {Bulletin of the American Mathematical Society}{Bull. Amer. Math. Soc.}}
   \ITEE{#2}{BAustrMS}{\ITE{\equal{#1}{+}}
      {Bulletin of the Australian Mathematical Society}{Bull. Austral. Math. Soc.}}
   \ITEE{#2}{BLondMS}{\ITE{\equal{#1}{+}}
      {Bulletin of the London Mathematical Sociecy}{Bull. Lond. Math. Soc.}}
   \ITEE{#2}{BAPolSSSM}{\ITE{\equal{#1}{+}}
      {Bulletin de l'Acad\'{e}mie Polonaise des Sciences. S\'{e}rie des Sciences Math\'{e}matiques}
      {Bull. Acad. Pol. Sci. S\'{e}r. Sci. Math.}}
   \ITEE{#2}{BullSM}{\ITE{\equal{#1}{+}}
      {Bulletin des Sciences Math\'{e}matiques}{Bull. Sci. Math.}}
   \ITEE{#2}{BullPol}{\ITE{\equal{#1}{+}}
      {Bulletin of the Polish Academy of Sciences: Mathematics}{Bull. Pol. Acad. Sci. Math.}}
   \ITEE{#2}{CanadJM}{\ITE{\equal{#1}{+}}
      {Canadian Journal Mathematics}{Canad. J. Math.}}
   \ITEE{#2}{CollectM}{\ITE{\equal{#1}{+}}
      {Collectanea Mathematica}{Collect. Math.}}
   \ITEE{#2}{CMUC}{\ITE{\equal{#1}{+}}
      {Commentationes Mathematicae Universitatis Carolinae}{Comment. Math. Univ. Carolin.}}
   \ITEE{#2}{CRParis}{\ITE{\equal{#1}{+}}
      {C. R. Paris}{C. R. Paris}}
   \ITEE{#2}{CRASParis}{\ITE{\equal{#1}{+}}
      {Comptes Rendus de l'Acad\'{e}mie des Sciences. Paris}{C. R. Acad. Sci. Paris}}
   \ITEE{#2}{CEurJM}{\ITE{\equal{#1}{+}}
      {Central European Journal of Mathematics}{Cent. Eur. J. Math.}}
   \ITEE{#2}{CMHelv}{\ITE{\equal{#1}{+}}
      {Commentarii Mathematici Helvetici}{Comment. Math. Helv.}}
   \ITEE{#2}{CollM}{\ITE{\equal{#1}{+}}
      {Colloquium Mathematicum}{Coll. Math.}}
   \ITEE{#2}{ComposM}{\ITE{\equal{#1}{+}}
      {Compositio Mathematica}{Compos. Math.}}
   \ITEE{#2}{CzMJ}{\ITE{\equal{#1}{+}}
      {Czechoslovak Mathematical Journal}{Czech. Math. J.}}
   \ITEE{#2}{DissM}{\ITE{\equal{#1}{+}}
      {Dissertationes Mathematicae}{Dissert. Math.}}
   \ITEE{#2}{DANSSSR}{\ITE{\equal{#1}{+}}
      {Doklady Akademii Nauk SSSR}{Dokl. Akad. Nauk SSSR}}
   \ITEE{#2}{DukeMJ}{\ITE{\equal{#1}{+}}
      {Duke Mathematical Journal}{Duke Math. J.}}
   \ITEE{#2}{ELA}{\ITE{\equal{#1}{+}}
      {The Electronic Journal of Linear Algebra}{Electron. J. Linear Algebra}}
   \ITEE{#2}{ExtrM}{\ITE{\equal{#1}{+}}
      {Extracta Mathematicae}{Extracta Math.}}
   \ITEE{#2}{FM}{\ITE{\equal{#1}{+}}
      {Fundamenta Mathematicae}{Fund. Math.}}
   \ITEE{#2}{FAA}{\ITE{\equal{#1}{+}}
      {Functional Analysis and its Applications}{Funct. Anal. Appl.}}
   \ITEE{#2}{FunkAnalPril}{\ITE{\equal{#1}{+}}
      {Funktsional'ny\u{\i} Analiz i Ego Prilozheniya}{Funkts. Anal. Prilozh.}}
   \ITEE{#2}{GTopA}{\ITE{\equal{#1}{+}}
      {General Topology and its Applications}{General Topol. Appl.}}
   \ITEE{#2}{HJM}{\ITE{\equal{#1}{+}}
      {Houston Journal of Mathematics}{Houston J. Math.}}
   \ITEE{#2}{IllinoisJM}{\ITE{\equal{#1}{+}}
      {Illinois Journal of Mathematics}{Illinois J. Math.}}
   \ITEE{#2}{IndagMP}{\ITE{\equal{#1}{+}}
      {Indagationes Mathematicae (Proceedings)}{Indagationes Math. Proc.}}
   \ITEE{#2}{IndianaUMJ}{\ITE{\equal{#1}{+}}
      {Indiana University Mathematical Journal}{Indiana Univ. Math. J.}}
   \ITEE{#2}{InHauEtSPM}{\ITE{\equal{#1}{+}}
      {Inst. Hautes \'{E}tudes Sci. Publ. Math.}{Inst. Hautes \'{E}tudes Sci. Publ. Math.}}
   \ITEE{#2}{IEOT}{\ITE{\equal{#1}{+}}
      {Integral Equations and Operator Theory}{Integr. Equ. Oper. Theory}}
   \ITEE{#2}{IsraelJM}{\ITE{\equal{#1}{+}}
      {Israel Journal of Mathematics}{Israel J. Math.}}
   \ITEE{#2}{JAusMSA}{\ITE{\equal{#1}{+}}
      {Journal of the Australian Mathematical Society. Series A}{J. Aust. Math. Soc. Ser. A}}
   \ITEE{#2}{JCA}{\ITE{\equal{#1}{+}}
      {Journal of Convex Analysis}{J. Convex Anal.}}
   \ITEE{#2}{JChinUST}{\ITE{\equal{#1}{+}}
      {J. China Univ. Sci. Tech.}{J. China Univ. Sci. Tech.}}
   \ITEE{#2}{JFA}{\ITE{\equal{#1}{+}}
      {Journal of Functional Analysis}{J. Funct. Anal.}}
   \ITEE{#2}{JKoreanMS}{\ITE{\equal{#1}{+}}
      {Journal of the Korean Mathematical Society}{J. Korean Math. Soc.}}
   \ITEE{#2}{JMAnApp}{\ITE{\equal{#1}{+}}
      {J. Math. Anal. Appl.}{J. Math. Anal. Appl.}}
   \ITEE{#2}{JOT}{\ITE{\equal{#1}{+}}
      {Journal of Operator Theory}{J. Operator Theory}}
   \ITEE{#2}{KodaiMSemRep}{\ITE{\equal{#1}{+}}
      {Kodai Math. Sem. Rep.}{Kodai Math. Sem. Rep.}}
   \ITEE{#2}{LAA}{\ITE{\equal{#1}{+}}
      {Linear Algebra and its Applications}{Linear Algebra Appl.}}
   \ITEE{#2}{LMLA}{\ITE{\equal{#1}{+}}
      {Linear and Multilinear Algebra}{Linear Multilinear Algebra}}
   \ITEE{#2}{LNM}{\ITE{\equal{#1}{+}}
      {Lecture Notes in Mathematics}{Lecture Notes Math.}}
   \ITEE{#2}{MathJap}{\ITE{\equal{#1}{+}}
      {Math. Japon.}{Math. Japon.}}
   \ITEE{#2}{MLQ}{\ITE{\equal{#1}{+}}
      {Mathematical Logic Quarterly}{Math. Log. Q.}}
   \ITEE{#2}{MProcCambPhS}{\ITE{\equal{#1}{+}}
      {Mathematical Proceedings of the Cambridge Philosophical Society}{Math. Proc. Cambridge Phil. Soc.}}
   \ITEE{#2}{MMag}{\ITE{\equal{#1}{+}}
      {Mathematics Magazine}{Math. Mag.}}
   \ITEE{#2}{MSb}{\ITE{\equal{#1}{+}}
      {Matematicheski\u{\i} Sbornik}{Mat. Sb.}}
   \ITEE{#2}{MStud}{\ITE{\equal{#1}{+}}
      {Matematychni Studi\"{\i}}{Mat. Stud.}}
   \ITEE{#2}{MScand}{\ITE{\equal{#1}{+}}
      {Mathematica Scandinavica}{Math. Scand.}}
   \ITEE{#2}{MAnn}{\ITE{\equal{#1}{+}}
      {Mathematische Annalen}{Math. Ann.}}
   \ITEE{#2}{MAMS}{\ITE{\equal{#1}{+}}
      {Memoirs of the American Mathematical Society}{Mem. Amer. Math. Soc.}}
   \ITEE{#2}{MichMJ}{\ITE{\equal{#1}{+}}
      {Michigan Mathematical Journal}{Mich. Math. J.}}
   \ITEE{#2}{MonatM}{\ITE{\equal{#1}{+}}
      {Monatshefte f\"{u}r Mathematik}{Mh. Math.}}
   \ITEE{#2}{NonlinA}{\ITE{\equal{#1}{+}}
      {Nonlinear Analysis: Theory, Methods \& Applications}{Nonlinear Anal.}}
   \ITEE{#2}{NAMS}{\ITE{\equal{#1}{+}}
      {Notices of the American Mathematical Society}{Notices Amer. Math. Soc.}}
   \ITEE{#2}{OpusM}{\ITE{\equal{#1}{+}}
      {Opuscula Mathematica}{Opuscula Math.}}
   \ITEE{#2}{PacJM}{\ITE{\equal{#1}{+}}
      {Pacific Journal of Mathematics}{Pacific J. Math.}}
   \ITEE{#2}{PeriodMHung}{\ITE{\equal{#1}{+}}
      {Periodica Mathematica Hungarica}{Period. Math. Hungarica}}
   \ITEE{#2}{PAMS}{\ITE{\equal{#1}{+}}
      {Proceedings of the American Mathematical Society}{Proc. Amer. Math. Soc.}}
   \ITEE{#2}{ProcCambPhS}{\ITE{\equal{#1}{+}}
      {Proceedings of the Cambridge Philosophical Society}{Proc. Cambridge Phil. Soc.}}
   \ITEE{#2}{ProcImpAcadTokyo}{\ITE{\equal{#1}{+}}
      {Proc. Imp. Acad. Tokyo}{Proc. Imp. Acad. Tokyo}}
   \ITEE{#2}{ProcKonink}{\ITE{\equal{#1}{+}}
      {Proceedings of the Koninklijke Nederlandse Akademie van Wetenschappen}{Nederl. Akad. Wetensch. Proc. Ser. A}}
   \ITEE{#2}{PLondMS}{\ITE{\equal{#1}{+}}
      {Proceedings of the London Mathematical Society}{Proc. London Math. Soc.}}
   \ITEE{#2}{PNatlUSA}{\ITE{\equal{#1}{+}}
      {Proceedings of the National Academy of Sciences of the United States of America}{Proc. Natl. Acad. Sci. USA}}
   \ITEE{#2}{PublRIMSKyoto}{\ITE{\equal{#1}{+}}
      {Publ. Res. Inst. Math. Sci. Kyoto Univ.}{Publ. Res. Inst. Math. Sci.}}
   \ITEE{#2}{PWN}{\ITE{\equal{#1}{+}}
      {PWN -- Polish Scientific Publishers, Warszawa}{PWN -- Polish Scientific Publishers, Warszawa}}
   \ITEE{#2}{RCMP}{\ITE{\equal{#1}{+}}
      {Rendiconti del Circolo Matematico di Palermo}{Rend. Circ. Mat. Palermo}}
   \ITEE{#2}{RussMS}{\ITE{\equal{#1}{+}}
      {Russian Mathematical Surveys}{Russian Math. Surveys}}
   \ITEE{#2}{SbM}{\ITE{\equal{#1}{+}}
      {Sbornik: Mathematics}{Sb. Math.}}
   \ITEE{#2}{SciRepTokyoA}{\ITE{\equal{#1}{+}}
      {Science Reports of Tokyo Kyoiku Daigaku, Section A}{Sci. Rep. Tokyo Kyoiku Daigaku Sect. A}}
   \ITEE{#2}{SeminProbStras}{\ITE{\equal{#1}{+}}
      {S\'{e}minaire de probabilit\'{e}s de Strasbourg}{S\'{e}min. Prob. Strasbourg}}
   \ITEE{#2}{SIAMJMAA}{\ITE{\equal{#1}{+}}
      {SIAM Journal on Matrix Analysis and Applications}{SIAM J. Matrix Anal. Appl.}}
   \ITEE{#2}{SibirMZ}{\ITE{\equal{#1}{+}}
      {Sibirski\v{\i} Mat. \v{Z}hurnal}{Sibirsk. Mat. \v{Z}.}}
   \ITEE{#2}{SM}{\ITE{\equal{#1}{+}}
      {Studia Mathematica}{Studia Math.}}
   \ITEE{#2}{TAMS}{\ITE{\equal{#1}{+}}
      {Transactions of the American Mathematical Society}{Trans. Amer. Math. Soc.}}
   \ITEE{#2}{TohokuMJ}{\ITE{\equal{#1}{+}}
      {T\^{o}hoku Mathematical Journal}{T\^{o}hoku Math. J.}}
   \ITEE{#2}{TomskUnivRev}{\ITE{\equal{#1}{+}}
      {Tomsk Universitet Review}{Tomsk. Univ. Rev.}}
   \ITEE{#2}{TopA}{\ITE{\equal{#1}{+}}
      {Topology and its Applications}{Topology Appl.}}
   \ITEE{#2}{TopMethNA}{\ITE{\equal{#1}{+}}
      {Topological Methods in Nonlinear Analysis}{Topol. Methods Nonlinear Anal.}}
   \ITEE{#2}{TsukubaJM}{\ITE{\equal{#1}{+}}
      {Tsukuba Journal of Mathematics}{Tsukuba J. Math.}}
   \ITEE{#2}{UspekhiMN}{\ITE{\equal{#1}{+}}
      {Uspekhi Matem. Nauk}{Uspekhi Mat. Nauk}}
   }

\newcommand{\paplist}[3][]{
   \ITEE{#3}{NIAkhiezer,IMGlazman1993}{
      \BIb{#2}{N.I. Akhiezer and I.M. Glazman}
         {Theory of Linear Operators in Hilbert Space}
         {Dover Publications, Inc., New York}{1993}{#1}}
   \ITEE{#3}{RDAnderson1967}{
      \BIB{#2}{R.D. Anderson}
         {On topological infinite deficiency}
         {\jRN{MichMJ}}{14}{1967}{365--383}{#1}}
   \ITEE{#3}{RDAnderson,JMcCharen1970}{
      \BIB{#2}{R.D. Anderson and J. McCharen}
         {On extending homeomorphisms to Fr\'{e}chet manifolds}
         {\jRN{PAMS}}{25}{1970}{283--289}{#1}}
   \ITEE{#3}{RDAnderson,DWCurtis,JVanMill1982}{
      \BIB{#2}{R.D. Anderson, D.W. Curtis, J. van Mill}
         {A fake topological Hilbert space}
         {\jRN{TAMS}}{272}{1982}{311--321}{#1}}
   \ITEE{#3}{RArens,JEells1956}{
      \BIB{#2}{R. Arens and J. Eells}
         {On embedding uniform and topological spaces}
         {\jRN{PacJM}}{6}{1956}{397--403}{#1}}
   \ITEE{#3}{NAronszajn,PPanitchpakdi1956}{
      \BIB{#2}{N. Aronszajn and P. Panitchpakdi}
         {Extension of uniformly continuous transformations and hyperconvex metric spaces}
         {\jRN{PacJM}}{6}{1956}{405--439}{#1}}
   \ITEE{#3}{KJBabenko1948}{
      \BIB{#2}{K.J. Babenko}
         {On conjugate functions}
         {\jRN{DANSSSR}}{62}{1948}{157--160}{#1}}
   \ITEE{#3}{TBanakh1995}{
      \BIB{#2}{T.O. Banakh}
         {Topology of spaces of probability measures, I}
         {\jRN{MStud}}{5}{1995}{65--87 (Russian)}{#1}}
   \ITEE{#3}{TBanakh1995a}{
      \BIB{#2}{T.O. Banakh}
         {Topology of spaces of probability measures, II}
         {\jRN{MStud}}{5}{1995}{88--106 (Russian)}{#1}}
   \ITEE{#3}{TBanakh1998}{
      \BIB{#2}{T. Banakh}
         {Characterization of spaces admitting a homotopy dense embedding into a Hilbert manifold}
         {\jRN{TopA}}{86}{1998}{123--131}{#1}}
   \ITEE{#3}{TBanakh,TNRadul1997}{
      \BIB{#2}{T.O. Banakh and T.N. Radul}
         {Topology of spaces of probability measures}
         {\jRN{SbM}}{188}{1997}{973--995}{#1}}
   \ITEE{#3}{TBanakh,TRadul,MZarichnyi1996}{
      \BIb{#2}{T. Banakh, T. Radul, M. Zarichnyi}
         {Absorbing sets in infinite-dimensional manifolds}
         {VNTL Publishers, Lviv}{1996}{#1}}
   \ITEE{#3}{TBanakh,IZarichnyy2008}{
      \BIB{#2}{T. Banakh and I. Zarichnyy}
         {Topological groups and convex sets homeomorphic to non-separable Hilbert spaces}
         {\jRN{CEurJM}}{6}{2008}{77--86}{#1}}
   \ITEE{#3}{HBecker,ASKechris1996}{
      \BIb{#2}{H. Becker and A.S. Kechris}
         {The Descriptive Set Theory of Polish Group Actions \textup{(London Math. Soc. Lecture Note Series, vol. 232)}}
         {University Press, Cambridge}{1996}{#1}}
   \ITEE{#3}{GBeer1993}{
      \BIb{#2}{G. Beer}
         {Topologies on Closed and Closed Convex Sets \textup{(Mathematics and Its Applications)}}
         {Kluwer Academic Publishers, Dordrecht}{1993}{#1}}
   \ITEE{#3}{NEBenamara,NNikolski1999}{
      \BIB{#2}{N.E. Benamara and N. Nikolski}
         {Resolvent tests for similarity to a normal operator}
         {\jRN{PLondMS}}{78}{1999}{585--626}{#1}}
   \ITEE{#3}{YBenyamini,JLindenstrauss2000}{
      \BIb{#2}{Y. Benyamini and J. Lindenstrauss}
         {Geometric nonlinear functional analysis I}
         {AMS Colloquium Publications 48}{2000}{#1}}
   \ITEE{#3}{SKBerberian1974}{
      \BIb{#2}{S.K. Berberian}
         {Lectures in Functional Analysis and Operator Theory}
         {Graduate Texts in Mathematics 15, Springer-Verlag, New York}{1974}{#1}}
   \ITEE{#3}{SNBernstein1954}{
      \BIb{#2}{S.N. Bernstein}
         {Collected Works II}
         {Akad. Nauk SSSR, Moscow}{1954 (Russian)}{#1}}
   \ITEE{#3}{CzBessaga,APelczynski1972}{
      \BIB{#2}{Cz. Bessaga and A. Pe\l{}czy\'{n}ski}
         {On spaces of measurable functions}
         {\jRN{SM}}{44}{1972}{597--615}{#1}}
   \ITEE{#3}{CzBessaga,APelczynski1975}{
      \BIb{#2}{Cz. Bessaga and A. Pe\l{}czy\'{n}ski}
         {Selected topics in infinite-dimensional topology}
         {\jRN{PWN}}{1975}{#1}}
   \ITEE{#3}{MBestvina,JMogilski1986}{
      \BIB{#2}{M. Bestvina and J. Mogilski}
         {Characterizing certain incomplete infinite-dimensional absolute retracts}
         {\jRN{MichMJ}}{33}{1986}{291--313}{#1}}
   \ITEE{#3}{MBestvina,PBowers,JMogilsky,JWalsh1986}{
      \BIB{#2}{M. Bestvina, P. Bowers, J. Mogilsky, J. Walsh}
         {Characterization of Hilbert space manifolds revisited}
         {\jRN{TopA}}{24}{1986}{53--69}{#1}}
   \ITEE{#3}{RBhatia1997}{
      \BIb{#2}{R. Bhatia}
         {Matrix Analysis}
         {Springer, New York}{1997}{#1}}
   \ITEE{#3}{GBirkhoff1936}{
      \BIB{#2}{G. Birkhoff}
         {A note on topological groups}
         {\jRN{ComposM}}{3}{1936}{427--430}{#1}}
   \ITEE{#3}{MSBirman,MZSolomjak1987}{
      \BIb{#2}{M.S. Birman and M.Z. Solomjak}
         {Spectral Theory of Self-Adjoint Operators in Hilbert Space}
         {D. Reidel Publishing Co., Dordrecht}{1987}{#1}}
   \ITEE{#3}{EBishop1961}{
      \BIB{#2}{E. Bishop}
         {A generalization of the Stone-Weierstrass theorem}
         {\jRN{PacJM}}{11}{1961}{777--783}{#1}}
   \ITEE{#3}{BBlackadar2006}{
      \BIb{#2}{B. Blackadar}{Operator Algebras. Theory of $\CCc^*$-algebras and von Neumann algebras 
         \textup{(Encyclopaedia of Mathematical Sciences, vol. 122: Operator Algebras and Non-Commutative Geometry III)}}
         {Springer-Verlag, Berlin-Heidelberg}{2006}{#1}}
   \ITEE{#3}{JBlass,WHolsztynski1972}{
      \BIB{#2}{J. Blass and W. Holszty\'{n}ski}
         {Cubical polyhedra and homotopy III}
         {\jRN{AttiAccLincRendNat}}{53}{1972}{275--279}{#1}}
   \ITEE{#3}{FFBonsall,NJDuncan1973}{
      \BIb{#2}{F.F. Bonsall and N.J. Duncan}
         {Complete Normed Algebras}
         {Springer Verlag, Berlin}{1973}{#1}}
   \ITEE{#3}{NBourbaki2002}{
      \BIb{#2}{N. Bourbaki}
         {Lie Groups and Lie Algebras, Chapters 4--6}
         {Springer, New York}{2002}{#1}}
   \ITEE{#3}{PLBowers1989}{
      \BIB{#2}{P.L. Bowers}
         {Limitation topologies on function spaces}
         {\jRN{TAMS}}{314}{1989}{421--431}{#1}}
   \ITEE{#3}{ABrown1953}{
      \BIB{#2}{A. Brown}
         {On a class of operators}
         {\jRN{PAMS}}{4}{1953}{723--728}{#1}}
   \ITEE{#3}{ABrown,CKFong,DWHadwin1978}{
      \BIB{#2}{A. Brown, C.-K. Fong, D.W. Hadwin}
         {Parts of operators on Hilbert space}
         {\jRN{IllinoisJM}}{22}{1978}{306--314}{#1}}
   \ITEE{#3}{AMBruckner,JBBruckner,BSThomson1997}{
      \BIb{#2}{A.M. Bruckner, J.B. Bruckner, B.S. Thomson}
         {Real Analysis}
         {Prentice-Hall, New Jersey}{1997}{#1}}
   \ITEE{#3}{PJCameron,AMVershik2006}{
      \BIB{#2}{P.J. Cameron and A.M. Vershik}
         {Some isometry groups of Urysohn space}
         {\jRN{AnnPALog}}{143}{2006}{70--78}{#1}}
   \ITEE{#3}{CCastaing1966}{
      \BIB{#2}{C. Castaing}
         {Quelques probl\`{e}mes de mesurabilit\'{e} li\'{e}es \`{a} la th\'{e}orie de la commande}
         {\jRN{CRParis}}{262}{1966}{409--411}{#1}}
   \ITEE{#3}{JAVanCasteren1980}{
      \BIB{#2}{J.A. van Casteren}
         {A problem of Sz.-Nagy}
         {\jRN{ActaSM}}{42}{1980}{189--194}{#1}}
   \ITEE{#3}{JAVanCasteren1983}{
      \BIB{#2}{J.A. van Casteren}
         {Operators similar to unitary or selfadjoint ones}
         {\jRN{PacJM}}{104}{1983}{241--255}{#1}}
   \ITEE{#3}{XCatepillan,MPtak,WSzymanski1994}{
      \BIB{#2}{X. Catepill\'{a}n, M. Ptak, W. Szyma\'{n}ski}
         {Multiple canonical decompositions of families of operators and a model of quasinormal families}
         {\jRN{PAMS}}{121}{1994}{1165--1172}{#1}}
   \ITEE{#3}{RCauty1994}{
      \BIB{#2}{R. Cauty}
         {Un espace m\'{e}trique lin\'{e}aire qui n'est pas un r\'{e}tracte absolu}
         {\jRN{FM}}{146}{1994}{85--99, (French)}{#1}}
   \ITEE{#3}{TAChapman1971}{
      \BIB{#2}{T.A. Chapman}
         {Deficiency in infinite-dimensional manifolds}
         {\jRN{GTopA}}{1}{1971}{263--272}{#1}}
   \ITEE{#3}{TAChapman1976}{
      \BIb{#2}{T.A. Chapman}
         {Lectures on Hilbert cube manifolds}
         {C.B.M.S. Regional Conference Series in Math. No 28, Amer. Math. Soc.}{1976}{#1}}
   \ITEE{#3}{RBChuaqui1977}{
      \BIB{#2}{R.B. Chuaqui}
         {Measures invariant under a group of transformations}
         {\jRN{PacJM}}{68}{1977}{313--329}{#1}}
   \ITEE{#3}{JBConway1985}{
      \BIb{#2}{J.B. Conway}
         {A Course in Functional Analysis}
         {Springer-Verlag, New York}{1985}{#1}}
   \ITEE{#3}{JBConway2000}{
      \BIb{#2}{J.B. Conway}
         {A Course in Operator Theory}
         {(Graduate Studies in Mathematics, vol. 21) Amer. Math. Soc., Providence}{2000}{#1}}
   \ITEE{#3}{GCorach,AMaestripieri,MMbekhta2009}{
      \BIB{#2}{G. Corach, A. Maestripieri, M. Mbekhta}
         {Metric and homogeneous structure of closed range operators}
         {\jRN{JOT}}{61}{2009}{171--190}{#1}}
   \ITEE{#3}{MJCowen,RGDouglas1978}{
      \BIB{#2}{M.J. Cowen and R.G. Douglas}
         {Complex geometry and operator theory}
         {\jRN{ActaM}}{141}{1978}{187--261}{#1}}
   \ITEE{#3}{DWCurtis1985}{
      \BIB{#2}{D.W. Curtis}
         {Boundary sets in the Hilbert cube}
         {\jRN{TopA}}{20}{1985}{201--221}{#1}}
   \ITEE{#3}{MMDay1958}{
      \BIb{#2}{M.M. Day}
         {Normed Linear Spaces}
         {Springer Verlag, Berlin}{1958}{#1}}
   \ITEE{#3}{CDellacherie1967}{
      \BIB{#2}{C. Dellacherie}
         {Un compl\'{e}ment au th\'{e}or\`{e}me de Weierstrass-Stone}
         {\jRN{SeminProbStras}}{1}{1967}{52--53}{#1}}
   \ITEE{#3}{JJDijkstra1987}{
      \BIB{#2}{J.J. Dijkstra}
         {Strong negligibility of $\sigma$-compacta does not characterize Hilbert space}
         {\jRN{PacJM}}{127}{1987}{19--30}{#1}}
   \ITEE{#3}{JJDijkstra1990}{
      \BIB{#2}{J.J. Dijkstra}
         {Characterizing Hilbert space topology in terms of strong negligibility}
         {\jRN{ComposM}}{75}{1990}{299--306}{#1}}
   \ITEE{#3}{TDobrowolski,WMarciszewski2002}{
      \BIB{#2}{T. Dobrowolski and W. Marciszewski}
         {Failure of the Factor Theorem for Borel pre-Hilbert spaces}
         {\jRN{FM}}{175}{2002}{53--68}{#1}}
   \ITEE{#3}{TDobrowolski,JMogilski1990}{
      \BiB{#2}{T. Dobrowolski and J. Mogilski}
         {Problems on Topological Classification of Incomplete Metric Spaces}{Chapter 25 in:}
         {Open Problems in Topology}{J. van Mill and G.M. Reed (eds.), North-Holland Amsterdam}{1990}{411--429}{#1}}
   \ITEE{#3}{TDobrowolski,HTorunczyk1981}{
      \BIB{#2}{T. Dobrowolski and H. Toru\'{n}czyk}
         {Separable complete ANR's admitting a group structure are Hilbert manifolds}
         {\jRN{TopA}}{12}{1981}{229--235}{#1}}
   \ITEE{#3}{RGDouglas1966}{
      \BIB{#2}{R.G. Douglas}
         {On majorization, factorization and range inclusion of operators in Hilbert space}
         {\jRN{PAMS}}{17}{1966}{413--416}{#1}}
   \ITEE{#3}{CHDowker1947}{
      \BIB{#2}{C.H. Dowker}
         {Mapping theorems for non-compact spaces}
         {\jRN{AmJM}}{69}{1947}{200--242}{#1}}
   \ITEE{#3}{CHDowker1952}{
      \BIB{#2}{C.H. Dowker}
         {Topology of metric complexes}
         {\jRN{AmJM}}{74}{1952}{555--577}{#1}}
   \ITEE{#3}{JDugundji1951}{
      \BIB{#2}{J. Dugundji}
         {An extension of Tietze's theorem}
         {\jRN{PacJM}}{1}{1951}{353--367}{#1}}
   \ITEE{#3}{JDugundji1958}{
      \BIB{#2}{J. Dugundji}
         {Absolute neighborhood retracts and local connectedness for arbitrary metric spaces}
         {\jRN{ComposM}}{13}{1958}{229--246}{#1}}
   \ITEE{#3}{JDugundji1965}{
      \BIB{#2}{J. Dugundji}
         {Locally equiconnected spaces and absolute neighborhood retracts}
         {\jRN{FM}}{57}{1965}{187--193}{#1}}
   \ITEE{#3}{NDunford,JTSchwartz1958}{
      \BIb{#2}{N. Dunford and J.T. Schwartz}
         {Linear Operators, part I}
         {Interscience Publishers, New York}{1958}{#1}}
   \ITEE{#3}{NDunford,JTSchwartz1963}{
      \BIb{#2}{N. Dunford and J.T. Schwartz}
         {Linear Operators, part II}
         {Interscience Publishers, New York}{1963}{#1}}
   \ITEE{#3}{NDunford,JTSchwartz1971}{
      \BIb{#2}{N. Dunford and J.T. Schwartz}
         {Linear Operators, part III}
         {Wiley-Interscience, New York}{1971}{#1}}
   \ITEE{#3}{MLEaton,MDPerlman1977}{
      \BIB{#2}{M.L. Eaton and M.D. Perlman}
         {Reflection groups, generalized Schur functions and the geometry of majorization}
         {\jRN{AnnProb}}{5}{1977}{829--860}{#1}}
   \ITEE{#3}{EGEffros1965}{
      \BIB{#2}{E.G. Effros}
         {The Borel space of von Neumann algebras on a separable Hilbert space}
         {\jRN{PacJM}}{15}{1965}{1153--1164}{#1}}
   \ITEE{#3}{EGEffros1966}{
      \BIB{#2}{E.G. Effros}
         {Global structure in von Neumann algebras}
         {\jRN{TAMS}}{121}{1966}{434--454}{#1}}
   \ITEE{#3}{REspinola,MAKhamsi2001}{
      \BiB{#2}{R. Espinola and M.A. Khamsi}
         {Introduction to hyperconvex spaces}{Chapter XIII in:}{Handbook of Metric Fixed Point Theory}
         {W.A. Kirk and B. Sims (editors), Kluwer Academic Publishers}{2001}{391--435}{#1}}
   \ITEE{#3}{PAFillmore,JPWilliams1971}{
      \BIB{#2}{P.A. Fillmore and J.P. Williams}
         {On operator ranges}
         {\jRN{AdvM}}{7}{1971}{254--281}{#1}}
   \ITEE{#3}{JEells,NHKuiper1969}{
      \BIB{#2}{J. Eells and N.H. Kuiper}
         {Homotopy negligible subsets in infinite-dimensional manifolds}
         {\jRN{ComposM}}{21}{1969}{151--161}{#1}}
   \ITEE{#3}{REngelking1977}{
      \BIb{#2}{R. Engelking}
         {General Topology}
         {\jRN{PWN}}{1977}{#1}}
   \ITEE{#3}{REngelking1978}{
      \BIb{#2}{R. Engelking}
         {Dimension Theory}
         {\jRN{PWN}}{1978}{#1}}
   \ITEE{#3}{REngelking1989}{
      \BIb{#2}{R. Engelking}
         {General Topology. Revised and completed edition \textup{(Sigma series in pure mathematics, vol. 6)}}
         {Heldermann Verlag, Berlin}{1989}{#1}}
   \ITEE{#3}{PErdos,RDMauldin1976}{
      \BIB{#2}{P. Erd\"{o}s and R.D. Mauldin}
         {The nonexistence of certain invariant measures}
         {\jRN{PAMS}}{59}{1976}{321--322}{#1}}
   \ITEE{#3}{JErnest1976}{
      \BIB{#2}{J. Ernest}
         {Charting the operator terrain}
         {\jRN{MAMS}}{171}{1976}{207 pp}{#1}}
   \ITEE{#3}{RHFox1943}{
      \BIB{#2}{R.H. Fox}
         {On fiber spaces, II}
         {\jRN{BAMS}}{49}{1943}{733--735}{#1}}
   \ITEE{#3}{NAFriedman1970}{
      \BIb{#2}{N.A. Friedman}
         {Introduction to ergodic theory}
         {Van Nostrand Reinhold Company}{1970}{#1}}
   \ITEE{#3}{MFujii,MKajiwara,YKato,FKubo1976}{
      \BIB{#2}{M. Fujii, M. Kajiwara, Y. Kato, F. Kubo}
         {Decompositions of operators in Hilbert spaces}
         {\jRN{MathJap}}{21}{1976}{117--120}{#1}}
   \ITEE{#3}{SGao,ASKechris2003}{
      \BIB{#2}{S. Gao and A.S. Kechris}
         {On the classification of Polish metric spaces up to isometry}
         {\jRN{MAMS}}{161}{2003}{viii+78}{#1}}
   \ITEE{#3}{MIGarrido,FMontalvo1991}{
      \BIB{#2}{M.I. Garrido and F. Montalvo}
         {On some generalizations of the Kakutani-Stone and Stone-Weierstrass theorems}
         {\jRN{ExtrM}}{6}{1991}{156--159}{#1}}
   \ITEE{#3}{LGe,JShen2002}{
      \BIB{#2}{L. Ge and J. Shen}
         {Generator problem for certain property T factors}
         {\jRN{PNAS}}{99}{2002}{565--567}{#1}}
   \ITEE{#3}{IMGelfand,MANaimark1943}{
      \BIB{#2}{I.M. Gelfand and M.A. Naimark}
         {On the embedding of normed rings into the ring of operators in Hilbert space}
         {\jRN{MSb}}{12}{1943}{197--213}{#1}}
   \ITEE{#3}{FGesztesy,MMalamud,MMitrea,SNaboko2009}{
      \BIB{#2}{F. Gesztesy, M. Malamud, M. Mitrea, S. Naboko}
         {Generalized polar decompositions for closed operators in Hilbert spaces and some applications}
         {\jRN{IEOT}}{64}{2009}{83--113}{#1}}
   \ITEE{#3}{LGillman,MJerison1960}{
      \BIb{#2}{L. Gillman and M. Jerison}
         {Rings of continuous functions}
         {New York}{1960}{#1}}
   \ITEE{#3}{JGlimm1960}{
      \BIB{#2}{J. Glimm}
         {A Stone-Weierstrass theorem for $\CCc^*$-algebras}
         {\jRN{AnnM}}{72}{1960}{216--244}{#1}}
   \ITEE{#3}{GGodefroy,NJKalton2003}{
      \BIB{#2}{G. Godefroy and N.J. Kalton}
         {Lipschitz-free Banach spaces}
         {\jRN{SM}}{159}{2003}{121--141}{#1}}
   \ITEE{#3}{ICGohberg,MGKrein1967}{
      \BIB{#2}{I.C. Gohberg and M.G. Krein}
         {On a description of contraction operators similar to unitary ones}
         {\jRN{FunkAnalPril}}{1}{1967}{38--60}{#1}}
   \ITEE{#3}{ELGriffinJr1953}{
      \BIB{#2}{E.L. Griffin Jr.}
         {Some contributions to the theory of rings of operators}
         {\jRN{TAMS}}{75}{1953}{471--504}{#1}}
   \ITEE{#3}{ELGriffinJr1955}{
      \BIB{#2}{E.L. Griffin Jr.}
         {Some contributions to the theory of rings of operators II}
         {\jRN{TAMS}}{79}{1955}{389--400}{#1}}
   \ITEE{#3}{MGromov1981}{
      \BIB{#2}{M. Gromov}
         {Groups of polynomial growth and expanding maps}
         {\jRN{InHauEtSPM}}{53}{1981}{53--73}{#1}}
   \ITEE{#3}{MGromov1999}{
      \BIb{#2}{M. Gromov}
         {Metric Structures for Riemannian and Non-Riemannian Spaces}
         {Progress in Math. \textbf{152}, Birkh\"{a}user}{1999}{#1}}
   \ITEE{#3}{JDeGroot1956}{
      \BIB{#2}{J. de Groot}
         {Non-archimedean metrics in topology}
         {\jRN{PAMS}}{7}{1956}{948--953}{#1}}
   \ITEE{#3}{LCGrove,CTBenson1985}{
      \BIb{#2}{L.C. Grove and C.T. Benson}
         {Finite Reflection Group}
         {2nd ed., Springer-Verlag}{1985}{#1}}
   \ITEE{#3}{VIGurarii1966}{
      \BIB{#2}{V.I. Gurari\v{\i}}{Spaces of universal placement, isotropic spaces and a problem of Mazur 
         on rotations of Banach spaces \textup{(Russian)}}
         {\jRN{SibirMZ}}{7}{1966}{1002--1013}{#1}}
   \ITEE{#3}{DWHadwin1976}{
      \BIB{#2}{D.W. Hadwin}
         {An operator-valued spectrum}
         {\jRN{NAMS}}{23}{1976}{A-163}{#1}}
   \ITEE{#3}{DWHadwin1977}{
      \BIB{#2}{D.W. Hadwin}
         {An operator-valued spectrum}
         {\jRN{IndianaUMJ}}{26}{1977}{329--340}{#1}}
   \ITEE{#3}{HHahn1932}{
      \BIb{#2}{H. Hahn}
         {Reelle Funktionen I}
         {Leipzig}{1932}{#1}}
   \ITEE{#3}{PRHalmos1950}{
      \BIb{#2}{P.R. Halmos}
         {Measure theory}
         {Van Nostrand, New York}{1950}{#1}}
   \ITEE{#3}{PRHalmos1951}{
      \BIb{#2}{P.R. Halmos}
         {Introduction to Hilbert Space and the Theory of Spectral Multiplicity}
         {Chelsea Publishing Company, New York}{1951}{#1}}
   \ITEE{#3}{PRHalmos1956}{
      \BIb{#2}{P.R. Halmos}
         {Lectures on Ergodic Theory}
         {Publ. Math. Soc. Japan, Tokyo}{1956}{#1}}
   \ITEE{#3}{PRHalmos1982}{
      \BIb{#2}{P.R. Halmos}
         {A Hilbert Space Problem Book}
         {Springer-Verlag New York Inc.}{1982}{#1}}
  \ITEE{#3}{PRHalmos,JEMcLaughlin1963}{
      \BIB{#2}{P.R. Halmos and J.E. McLaughlin}
         {Partial isometries}
         {\jRN{PacJM}}{13}{1963}{585--596}{#1}}
   \ITEE{#3}{RWHansell1972}{
      \BIB{#2}{R.W. Hansell}
         {On the nonseparable theory of Borel and Souslin sets}
         {\jRN{BAMS}}{78}{1972}{236--241}{#1}}
   \ITEE{#3}{FHausdorff1930}{
      \BIB{#2}{F. Hausdorff}
         {Erweiterung einer Hom\"{o}omorphie}
         {\jRN{FM}}{16}{1930}{353--360}{#1}}
   \ITEE{#3}{FHausdorff1934}{
      \BIB{#2}{F. Hausdorff}
         {\"{U}ber innere Abbildungen}
         {\jRN{FM}}{23}{1934}{279--291}{#1}}
   \ITEE{#3}{FHausdorff1938}{
      \BIB{#2}{F. Hausdorff}
         {Erweiterung einer stetigen Abbildung}
         {\jRN{FM}}{30}{1938}{40--47}{#1}}
   \ITEE{#3}{DWHenderson1971}{
      \BIB{#2}{D.W. Henderson}
         {Corrections and extensions of two papers about infinite-dimensional manifolds}
         {\jRN{GTopA}}{1}{1971}{321--327}{#1}}
   \ITEE{#3}{DWHenderson1975}{
      \BIB{#2}{D.W. Henderson}
         {$Z$-sets in ANR's}
         {\jRN{TAMS}}{213}{1975}{205--216}{#1}}
   \ITEE{#3}{DWHenderson,RMSchori1970}{
      \BIB{#2}{D.W. Henderson and R.M. Schori}
         {Topological classification of infinite-dimensional manifolds by homotopy type}
         {\jRN{BAMS}}{76}{1970}{121--124}{#1}}
   \ITEE{#3}{DWHenderson,JEWest1970}{
      \BIB{#2}{D.W. Henderson and J.E. West}
         {Triangulated infinite-dimensional manifolds}
         {\jRN{BAMS}}{76}{1970}{655--660}{#1}}
   \ITEE{#3}{BHoffmann1979}{
      \BIB{#2}{B. Hoffmann}
         {A compact contractible topological group is trivial}
         {\jRN{ArchM}}{32}{1979}{585--587}{#1}}
   \ITEE{#3}{DHofmann2002}{
      \BIB{#2}{D. Hofmann}
         {On a generalization of the Stone-Weierstrass theorem}
         {\jRN{ACS}}{10}{2002}{569--592}{#1}}
   \ITEE{#3}{GHognas,AMukherjea1995}{
      \BIb{#2}{G. H\"ogn\"as and A. Mukherjea}
         {Probability Measures on Semigroups. Convolution Products, Random Walks, and Random Matrices}
         {Plenum Press, New York}{1995}{#1}}
   \ITEE{#3}{MRHolmes1992}{
      \BIB{#2}{M.R. Holmes}
         {The universal separable metric space of Urysohn and isometric embeddings thereof in Banach spaces}
         {\jRN{FM}}{140}{1992}{199--223}{#1}}
   \ITEE{#3}{MRHolmes2008}{
      \BIB{#2}{M.R. Holmes}
         {The Urysohn space embeds in Banach spaces in just one way}
         {\jRN{TopA}}{155}{2008}{1479--1482}{#1}}
   \ITEE{#3}{RRHolmes,TYTam1999}{
      \BIB{#2}{R.R. Holmes and T.Y. Tam}
         {Distance to the convex hull of an orbit under the action of a compact group}
         {\jRN{JAusMSA}}{66}{1999}{331--357}{#1}}
   \ITEE{#3}{RHorn,RMathias1990}{
      \BIB{#2}{R. Horn and R. Mathias}
         {Cauchy-Schwartz inequalities associated with positive semidefinite matrices}
         {\jRN{LAA}}{142}{1990}{63--82}{#1}}
   \ITEE{#3}{GEHuhunaisvili1955}{
      \BIB{#2}{G.E. Huhunai\v{s}vili}
         {On a property of Urysohn's universal metric space}
         {\jRN{DANSSSR}}{101}{1955}{607--610 (Russian)}{#1}}
   \ITEE{#3}{JEHumphreys1990}{
      \BIb{#2}{J.E. Humphreys}
         {Reflection Groups and Coxeter Groups}
         {Cambridge University Press}{1990}{#1}}
   \ITEE{#3}{JRIsbell1964}{
      \BIB{#2}{J.R. Isbell}
         {Six theorems about injective metric spaces}
         {\jRN{CMHelv}}{39}{1964}{65--76}{#1}}
   \ITEE{#3}{SIzumino,YKato1985}{
      \BIB{#2}{S. Izumino and Y. Kato}
         {The closure of invertible operators on Hilbert space}
         {\jRN{ActaSM}}{49}{1985}{321--327}{#1}}
   \ITEE{#3}{CJiang2004}{
      \BIB{#2}{C. Jiang}
         {Similarity classification of Cowen-Douglas operators}
         {\jRN{CanadJM}}{56}{2004}{742--775}{#1}}
   \ITEE{#3}{WBJohnson,JLindenstrauss2001}{
      \BiB{#2}{W.B. Johnson and J. Lindenstrauss}{Basic Concepts in the Geometry of Banach Spaces}
         {Chapter 1 in:}{Handbook of the Geometry of Banach Spaces, Vol. 1}
         {W.B. Johnson and J. Lindenstrauss (editors), Elsevier Science B.V., Amsterdam}{2001}{1--84}{#1}}
   \ITEE{#3}{IBJung,JStochel2008}{
      \BIB{#2}{I.B. Jung and J. Stochel}
         {Subnormal operators whose adjoints have rich point spectrum}
         {\jRN{JFA}}{255}{2008}{1797--1816}{#1}}
   \ITEE{#3}{RVKadison,JRRingrose1983}{
      \BIb{#2}{R.V. Kadison and J.R. Ringrose}
         {Fundamentals of the Theory of Operator Algebras. Volume I: Elementary Theory}
         {Academic Press, Inc., New York-London}{1983}{#1}}
   \ITEE{#3}{RVKadison,JRRingrose1986}{
      \BIb{#2}{R.V. Kadison and J.R. Ringrose}
         {Fundamentals of the Theory of Operator Algebras. Volume II: Advanced Theory}
         {Academic Press, Inc., Orlando-London}{1986}{#1}}
   \ITEE{#3}{SKakutani1936}{
      \BIB{#2}{S. Kakutani}
         {\"{U}ber die Metrisation der topologischen Gruppen}
         {\jRN{ProcImpAcadTokyo}}{12}{1936}{82--84}{#1}}
   \ITEE{#3}{SKakutani1938}{
      \BIB{#2}{S. Kakutani}
         {Two fixed-point theorems concerning bicompact convex sets}
         {\jRN{ProcImpAcadTokyo}}{14}{1938}{242--245}{#1}}
   \ITEE{#3}{SKakutani1941}{
      \BIB{#2}{S. Kakutani}
         {Concrete representation of abstract L-spaces}
         {\jRN{AnnM}}{42}{1941}{523--537}{#1}}
   \ITEE{#3}{SKakutani1941a}{
      \BIB{#2}{S. Kakutani}
         {Concrete representation of abstract M-spaces}
         {\jRN{AnnM}}{42}{1941}{994--1024}{#1}}
   \ITEE{#3}{NKalton2007}{
      \BIB{#2}{N. Kalton}
         {Extending Lipschitz maps into $\CCc(K)$-spaces}
         {\jRN{IsraelJM}}{162}{2007}{275--315}{#1}}
   \ITEE{#3}{RKane2001}{
      \BIb{#2}{R. Kane}
         {Reflection Groups and Invariant Theory}
         {Canadian Mathematical Society, Springer}{2001}{#1}}
   \ITEE{#3}{VKannan,SRRaju1980}{
      \BIB{#2}{V. Kannan and S.R. Raju}
         {The nonexistence of invariant universal measures on semigroups}
         {\jRN{PAMS}}{78}{1980}{482--484}{#1}}
   \ITEE{#3}{IKaplansky1951}{
      \BIB{#2}{I. Kaplansky}
         {A theorem on rings of operators}
         {\jRN{PacJM}}{1}{1951}{227--232}{#1}}
   \ITEE{#3}{MKatetov1988}{
      \BiB{#2}{M. Kat\v{e}tov}{On universal metric spaces}{in: Frolik (ed.),}
         {General Topology and its Relations to Modern Analysis and Algebra VI. Proceedings of the Sixth Prague 
         Topological Symposium 1986}{Heldermann Verlag Berlin}{1988}{323--330}{#1}}
   \ITEE{#3}{YKatznelson1960}{
      \BIB{#2}{Y. Katznelson}
         {Sur les alg\'{e}bres dont les \'{e}l\'{e}ments non n\'{e}gatifs admettent des racines carr\'{e}es}
         {\jRN{AnnSciEcNormSupT}}{77}{1960}{167--174}{#1}}
   \ITEE{#3}{OHKeller1931}{
      \BIB{#2}{O.H. Keller}
         {Die Homoiomorphie der kompakten konvexen Mengen in Hilbertschen Raum}
         {\jRN{MAnn}}{105}{1931}{748--758}{#1}}
   \ITEE{#3}{MAKhamsi,WAKirk,CMartinez2000}{
      \BIB{#2}{M.A. Khamsi, W.A. Kirk, C. Martinez}
         {Fixed point and selection theorems in hyperconvex spaces}
         {\jRN{PAMS}}{128}{2000}{3275--3283}{#1}}
   \ITEE{#3}{ABKhararazishvili1998}{
      \BIb{#2}{A.B. Khararazishvili}
         {Transformation groups and invariant measures. Set-theoretic aspects}
         {World Scientific Publishing Co., Inc., River Edge, NJ}{1998}{#1}}
   \ITEE{#3}{YKijima1987}{
      \BIB{#2}{Y. Kijima}
         {Fixed points of nonexpansive self-maps of a compact metric space}
         {\jRN{JMAnApp}}{123}{1987}{114--116}{#1}}
  \ITEE{#3}{JSKim,ChRKim,SGLee1980}{
      \BIB{#2}{J.S. Kim, Ch.R. Kim, S.G. Lee}
         {Reducing operator valued spectra of a Hilbert space operator}
         {\jRN{JKoreanMS}}{17}{1980}{123--129}{#1}}
   \ITEE{#3}{JKindler1995}{
      \BIB{#2}{J. Kindler}
         {Minimax theorems with applications to convex metric spaces}
         {\jRN{CollM}}{68}{1995}{179--186}{#1}}
   \ITEE{#3}{WAKirk1998}{
      \BIB{#2}{W.A. Kirk}
         {Hyperconvexity of $\RRR$-trees}
         {\jRN{FM}}{156}{1998}{67--72}{#1}}
   \ITEE{#3}{VLKleeJr1952}{
      \BIB{#2}{V.L. Klee Jr.}
         {Invariant metrics in groups (solution of a problem of Banach)}
         {\jRN{PAMS}}{3}{1952}{484--487}{#1}}
   \ITEE{#3}{HJKowalsky1957}{
      \BIB{#2}{H.J. Kowalsky}
         {Einbettung metrischer R\"{a}ume}
         {\jRN{ArchM}}{8}{1957}{336--339}{#1}}
   \ITEE{#3}{WKubis,MRubin2010}{
      \BIB{#2}{W. Kubi\'{s} and M. Rubin}
         {Extension and reconstruction theorems for the Urysohn universal metric space}
         {\jRN{CzMJ}}{60}{2010}{1--29}{#1}}
   \ITEE{#3}{KKuratowski1966}{
      \BIb{#2}{K. Kuratowski}
         {Topology. \textup{Vol. I}}
         {\jRN{PWN}}{1966}{#1}}
   \ITEE{#3}{KKuratowski,BKnaster1927}{
      \BIB{#2}{K. Kuratowski and B. Knaster}
         {A connected and connected im kleinen point set which contains no perfect subset}
         {\jRN{BAMS}}{33}{1927}{106--109}{#1}}
   \ITEE{#3}{KKuratowski,AMostowski1976}{
      \BIb{#2}{K. Kuratowski and A. Mostowski}
         {Set Theory with an Introduction to Descriptive Set Theory}
         {\jRN{PWN}}{1976}{#1}}
   \ITEE{#3}{GLewicki1992}{
      \BIB{#2}{G. Lewicki}
         {Bernstein's ``lethargy'' theorem in metrizable topological linear spaces}
         {\jRN{MonatM}}{113}{1992}{213--226}{#1}}
   \ITEE{#3}{ASLewis1996}{
      \BIB{#2}{A.S. Lewis}
         {Group invariance and convex matrix analysis}
         {\jRN{SIAMJMAA}}{17}{1996}{927--949}{#1}}
   \ITEE{#3}{C-KLi,N-KTsing1991}{
      \BIB{#2}{C.-K. Li and N.-K. Tsing}
         {$G$-invariant norms and $G(c)$-radii}
         {\jRN{LAA}}{150}{1991}{179--194}{#1}}
   \ITEE{#3}{AJLazar,JLindenstrauss1971}{
      \BIB{#2}{A.J. Lazar and J. Lindenstrauss}
         {Banach spaces whose duals are $L_1$ spaces and their representing matrices}
         {\jRN{ActaM}}{126}{1971}{165--193}{#1}}
   \ITEE{#3}{EHLieb,MLoss1997}{
      \BIb{#2}{E.H. Lieb and M. Loss}
         {Analysis \textup{(Graduate Studies in Mathematics, vol. 14)}}
         {Amer. Math. Soc., Providence, RI}{1997}{#1}}
   \ITEE{#3}{ALindenbaum1926}{
      \BIB{#2}{A. Lindenbaum}
         {Contributions \`{a} l'\'{e}tude de l'espace m\'{e}trique I}
         {\jRN{FM}}{8}{1926}{209--222}{#1}}
   \ITEE{#3}{DLindenstrauss,LTzafriri1971}{
      \BIB{#2}{D. Lindenstrauss and L. Tzafriri}
         {On the complemented subspaces problem}
         {\jRN{IsraelJM}}{9}{1971}{263--269}{#1}}
   \ITEE{#3}{RILoebl1986}{
      \BIB{#2}{R.I. Loebl}
         {A note on containment of operators}
         {\jRN{BAustrMS}}{33}{1986}{279--291}{#1}}
   \ITEE{#3}{LHLoomis1945}{
      \BIB{#2}{L.H. Loomis}
         {Abstract congruence and the uniqueness of Haar measure}
         {\jRN{AnnM}}{46}{1945}{348--355}{#1}}
   \ITEE{#3}{LHLoomis1949}{
      \BIB{#2}{L.H. Loomis}
         {Haar measure in uniform structures}
         {\jRN{DukeMJ}}{16}{1949}{193--208}{#1}}
   \ITEE{#3}{ERLorch1939}{
      \BIB{#2}{E.R. Lorch}
         {Bicontinuous linear transformation in certain vector spaces}
         {\jRN{BAMS}}{45}{1939}{564--569}{#1}}
   \ITEE{#3}{ATLundell,SWeingram1969}{
      \BIb{#2}{A.T. Lundell and S. Weingram}
         {The topology of CW-complexes}
         {Litton Educ. Publ.}{1969}{#1}}
   \ITEE{#3}{WLusky1976}{
      \BIB{#2}{W. Lusky}
         {The Gurarij spaces are unique}
         {\jRN{ArchM}}{27}{1976}{627--635}{#1}}
   \ITEE{#3}{WLusky1977}{
      \BIB{#2}{W. Lusky}
         {On separable Lindenstrauss spaces}
         {\jRN{JFA}}{26}{1977}{103--120}{#1}}
   \ITEE{#3}{DMaharam1942}{
      \BIB{#2}{D. Maharam}
         {On homogeneous measure algebras}
         {\jRN{PNatlUSA}}{28}{1942}{108--111}{#1}}
   \ITEE{#3}{MMalicki,SSolecki2009}{
      \BIB{#2}{M. Malicki and S. Solecki}
         {Isometry groups of separable metric spaces}
         {\jRN{MProcCambPhS}}{146}{2009}{67--81}{#1}}
   \ITEE{#3}{PMankiewicz1972}{
      \BIB{#2}{P. Mankiewicz}
         {On extension of isometries in normed linear spaces}
         {\jRN{BAPolSSSM}}{20}{1972}{367--371}{#1}}
   \ITEE{#3}{JMartinezMaurica,MTPellon1987}{
      \BIB{#2}{J. Martinez-Maurica and M.T. Pell\'{o}n}
         {Non-archimedean Chebyshev centers}
         {\jRN{IndagMP}}{90}{1987}{417--421}{#1}}
   \ITEE{#3}{KMaurin1980}{
      \BIb{#2}{K. Maurin}
         {Analysis, Part II}
         {D. Reidel, Dordrecht-Boston-London}{1980}{#1}}
   \ITEE{#3}{SMazur,SUlam1932}{
      \BIB{#2}{S. Mazur and S. Ulam}
         {Sur les transformationes isom\'{e}triques d'espaces vectoriels norm\'{e}s}
         {\jRN{CRASParis}}{194}{1932}{946--948}{#1}}
   \ITEE{#3}{SMazurkiewicz1920}{
      \BIB{#2}{S. Mazurkiewicz}
         {Sur les lignes de Jordan}
         {\jRN{FM}}{1}{1920}{166--209}{#1}}
   \ITEE{#3}{SMazurkiewicz,WSierpinski1920}{
      \BIB{#2}{S. Mazurkiewicz and W. Sierpi\'{n}ski}
         {Contributions a la topologie des ensembles denombrables}
         {\jRN{FM}}{1}{1920}{17--27}{#1}}
   \ITEE{#3}{MMbekhta1992}{
      \BIB{#2}{M. Mbekhta}
         {Sur la structure des composantes connexes semi-Fredholm de $B(H)$}
         {\jRN{PAMS}}{116}{1992}{521--524}{#1}}
   \ITEE{#3}{JEMcCarthy1996}{
      \BIB{#2}{J.E. McCarthy}
         {Boundary values and Cowen-Douglas curvature}
         {\jRN{JFA}}{137}{1996}{1--18}{#1}}
   \ITEE{#3}{JMelleray2007}{
      \BIB{#2}{J. Melleray}
         {Computing the complexity of the relation of isometry between separable Banach spaces}
         {\jRN{MLQ}}{53}{2007}{128--131}{#1}}
   \ITEE{#3}{JMelleray2007a}{
      \BIB{#2}{J. Melleray}
         {On the geometry of Urysohn's universal metric space}
         {\jRN{TopA}}{154}{2007}{384--403}{#1}}
   \ITEE{#3}{JMelleray2008}{
      \BIB{#2}{J. Melleray}
         {Some geometric and dynamical properties of the Urysohn space}
         {\jRN{TopA}}{155}{2008}{1531--1560}{#1}}
   \ITEE{#3}{JMelleray,FVPetrov,AMVershik2008}{
      \BIB{#2}{J. Melleray, F.V. Petrov, A.M. Vershik}
         {Linearly rigid metric spaces and the embedding problem}
         {\jRN{FM}}{199}{2008}{177--194}{#1}}
   \ITEE{#3}{EMichael1953}{
      \BIB{#2}{E. Michael}
         {Some extension theorems for continuous functions}
         {\jRN{PacJM}}{3}{1953}{789--806}{#1}}
   \ITEE{#3}{EMichael1954}{
      \BIB{#2}{E. Michael}
         {Local properties of topological spaces}
         {\jRN{DukeMJ}}{21}{1954}{163--171}{#1}}
   \ITEE{#3}{EMichael1956}{
      \BIB{#2}{E. Michael}
         {Selected selection theorems}
         {\jRN{AmMMon}}{58}{1956}{233--238}{#1}}
   \ITEE{#3}{EMichael1956a}{
      \BIB{#2}{E. Michael}
         {Continuous selections. I}
         {\jRN{AnnM}}{63}{1956}{361--382}{#1}}
   \ITEE{#3}{EMichael1956b}{
      \BIB{#2}{E. Michael}
         {Continuous selections. II}
         {\jRN{AnnM}}{64}{1956}{562--580}{#1}}
   \ITEE{#3}{EMichael1959}{
      \BIB{#2}{E. Michael}
         {A theorem on semi-continuous set-valued functions}
         {\jRN{DukeMJ}}{26}{1959}{647--652}{#1}}
   \ITEE{#3}{JVanMill1986}{
      \BIB{#2}{J. van Mill}
         {Another counterexample in ANR theory}
         {\jRN{PAMS}}{97}{1986}{136--138}{#1}}
   \ITEE{#3}{JVanMill2001}{
      \BIb{#2}{J. van Mill}
         {The Infinite-Dimensional Topology of Function Spaces 
         \textup{(North-Holland Mathematical Library, vol. 64)}}
         {Elsevier, Amsterdam}{2001}{#1}}
   \ITEE{#3}{WMlak1991}{
      \BIb{#2}{W. Mlak}
         {Hilbert Spaces and Operator Theory}
         {PWN --- Polish Scientific Publishers and Kluwer Academic Publishers, Warszawa-Dordrecht}{1991}{#1}}
   \ITEE{#3}{JMogilski1979}{
      \BIB{#2}{J. Mogilski}
         {$CE$-decomposition of $l_2$-manifolds}
         {\jRN{BAPolSSSM}}{27}{1979}{309--314}{#1}}
   \ITEE{#3}{RLMoore1916}{
      \BIB{#2}{R.L. Moore}
         {On the foundations of plane analysis situs}
         {\jRN{TAMS}}{17}{1916}{131--164}{#1}}
   \ITEE{#3}{KMorita1955}{
      \BIB{#2}{K. Morita}
         {A condition for the metrizability of topological spaces and for $n$-dimensionality}
         {\jRN{SciRepTokyoA}}{5}{1955}{33--36}{#1}}
   \ITEE{#3}{AMukherjea,NATserpes1976}{
      \BIb{#2}{A. Mukherjea and N.A. Tserpes}
         {Measures on topological semigroups}
         {Springer Lecture Notes in Math. Vol. 547, Berlin}{1976}{#1}}
   \ITEE{#3}{JMycielski1974}{
      \BIB{#2}{J. Mycielski}
         {Remarks on invariant measures in metric spaces}
         {\jRN{CollM}}{32}{1974}{105--112}{#1}}
   \ITEE{#3}{SNNaboko1984}{
      \BIB{#2}{S.N. Naboko}
         {Conditions for similarity to unitary and selfadjoint operators}
         {\jRN{FunkAnalPril}}{18}{1984}{16--27}{#1}}
   \ITEE{#3}{LNachbin1965}{
      \BIb{#2}{L. Nachbin}
         {The Haar Integral}
         {D. Van Nostrand Company, Inc., Princeton-New Jersey-Toronto-New York-London}{1965}{#1}}
   \ITEE{#3}{TDNarang,SKGarg1991}{
      \BIB{#2}{T.D. Narang and S.K. Garg}
         {On the uniqueness of best approximation in non-archimedian spaces}
         {\jRN{PeriodMHung}}{22}{1991}{121--124}{#1}}
   \ITEE{#3}{JVonNeumann1930}{
      \BIB{#2}{J. von Neumann}
         {Zur Algebra der Funktionaloperationen und Theorie der normalen Operatoren}
         {\jRN{MAnn}}{102}{1930}{370--427}{#1}}
   \ITEE{#3}{JVonNeumann1934}{
      \BIB{#2}{J. von Neumann}
         {Zum Haarschen Mass in topologischen Gruppen}
         {\jRN{ComposM}}{1}{1934}{106--114}{#1}}
   \ITEE{#3}{JVonNeumann1937}{
      \BiB{#2}{J. von Neumann}
         {Some matrix-inequalities and metrization of matrix-space}{\jRN{TomskUnivRev}{} \textbf{1} (1937), 286--300; 
         in }{Collected Works}{Pergamon, New York}{1962}{Vol. 4, 205--219}{#1}}
   \ITEE{#3}{JVonNeumann1949}{
      \BIB{#2}{J. von Neumann}
         {On Rings of Operators. Reduction Theory}
         {\jRN{AnnM}}{50}{1949}{401--485}{#1}}
   \ITEE{#3}{ONielson1973}{
      \BIB{#2}{O. Nielson}
         {Borel sets of von Neumann algebras}
         {\jRN{AmJM}}{95}{1973}{145--164}{#1}}
   \ITEE{#3}{pn2002}{\bibITEM{#2}{#1} \mypaplist{pn1}}
   \ITEE{#3}{pn2006a}{\bibITEM{#2}{#1} \mypaplist{pn2}}
   \ITEE{#3}{pn2006b}{\bibITEM{#2}{#1} \mypaplist{pn3}}
   \ITEE{#3}{pn2007}{\bibITEM{#2}{#1} \mypaplist{pn4}}
   \ITEE{#3}{pn2008a}{\bibITEM{#2}{#1} \mypaplist{pn5}}
   \ITEE{#3}{pn2008b}{\bibITEM{#2}{#1} \mypaplist{pn6}}
   \ITEE{#3}{pn2009a}{\bibITEM{#2}{#1} \mypaplist{pn7}}
   \ITEE{#3}{pn2009b}{\bibITEM{#2}{#1} \mypaplist{pn8}}
   \ITEE{#3}{pn2009c}{\bibITEM{#2}{#1} \mypaplist{pn9}}
   \ITEE{#3}{pn2010a}{\bibITEM{#2}{#1} \mypaplist{pn12}}
   \ITEE{#3}{pn2010b}{\bibITEM{#2}{#1} \mypaplist{pn13}}
   \ITEE{#3}{pn2011a}{\bibITEM{#2}{#1} \mypaplist{pn10}}
   \ITEE{#3}{pn2011b}{\bibITEM{#2}{#1} \mypaplist{pn15}}
   \ITEE{#3}{pn2011c}{\bibITEM{#2}{#1} \mypaplist{pn16}}
   \ITEE{#3}{pn2011d}{\bibITEM{#2}{#1} \mypaplist{pn17}}
   \ITEE{#3}{pn2009x}{
      \bibITEM{#2}{#1} \mypaplist{pn11}}
   \ITEE{#3}{pn2010x}{
      \bibITEM{#2}{#1} \mypaplist{pn14}}
   \ITEE{#3}{pnXXXXb}{
      \bibITEM{#2}{#1} \mypaplist{pnX2}}
   \ITEE{#3}{pnXXXXc}{
      \bibITEM{#2}{#1} \mypaplist{pnX3}}
   \ITEE{#3}{pnXXXXd}{
      \bibITEM{#2}{#1} \mypaplist{pnX13}}
   \ITEE{#3}{MNiezgoda1998}{
      \BIB{#2}{M. Niezgoda}
         {Group majorization and Schur type inequalities}
         {\jRN{LAA}}{268}{1998}{9--30}{#1}}
   \ITEE{#3}{MNiezgoda1998a}{
      \BIB{#2}{M. Niezgoda}
         {An analytical characterization of effective and of irreducible groups inducing cone orderings}
         {\jRN{LAA}}{269}{1998}{105--114}{#1}}
   \ITEE{#3}{MNiezgoda,TYTam2001}{
      \BIB{#2}{M. Niezgoda and T.Y. Tam}
         {On norm property of $G(c)$-radii and Eaton triples}
         {\jRN{LAA}}{336}{2001}{119--130}{#1}}
   \ITEE{#3}{APazy1983}{
      \BIb{#2}{A. Pazy}{Semigroups of Linear Operators 
         and Applications to Partial Differential Equations \textup{(Applied Mathematical Sciences, vol. 44)}}
         {Springer-Verlag, New York}{1983}{#1}}
   \ITEE{#3}{APelc1982}{
      \BIB{#2}{A. Pelc}
         {Semiregular invariant measures on abelian groups}
         {\jRN{PAMS}}{86}{1982}{423--426}{#1}}
   \ITEE{#3}{RPenrose1955}{
      \BIB{#2}{R. Penrose}
         {A generalized inverse for matrices}
         {\jRN{ProcCambPhS}}{51}{1955}{406--413}{#1}}
   \ITEE{#3}{VPestov2006}{
      \BIb{#2}{V. Pestov}
         {Dynamics of infinite-dimensional groups. The Ramsey-Dvoretzky-Milman phenomenon}
         {University Lecture Series \textbf{40}, AMS, Providence, RI}{2006}{#1}}
   \ITEE{#3}{VPestov2007}{
      \BiB{#2}{V. Pestov}
         {Forty-plus annotated questions about large topological groups}
         {in:}{Open Problems in Topology II}{Elliot Pearl (editor), Elsevier B.V., Amsterdam}{2007}{439--450}{#1}}
   \ITEE{#3}{PVPetersen1993}{
      \BiB{#2}{P.V. Petersen}
         {Gromov-Hausdorff convergence of metric spaces}{in book:}{Differential Geometry: Riemannian Geometry 
         (Los Angeles, CA, 1990)}{Amer. Math. Soc., Providence, RI}{1993}{489--504}{#1}}
   \ITEE{#3}{DRamachandran,MMisiurewicz1982}{
      \BIB{#2}{D. Ramachandran and M. Misiurewicz}
         {Hopf's theorem on invariant measures for a group of transformations}
         {\jRN{SM}}{74}{1982}{183--189}{#1}}
   \ITEE{#3}{JMRosenblatt1974}{
      \BIB{#2}{J.M. Rosenblatt}
         {Equivalent invariant measures}
         {\jRN{IsraelJM}}{17}{1974}{261--270}{#1}}
   \ITEE{#3}{HLRoyden1963}{
      \BIb{#2}{H.L. Royden}
         {Real Analysis}
         {The Macmillan Co., New York}{1963}{#1}}
   \ITEE{#3}{WRudin1962}{
      \BIb{#2}{W. Rudin}
         {Fourier Analysis on Groups \textup{(Interscience Tracts in Pure and Applied Mathematics, Number 12)}}
         {Interscience Publishers, New York}{1962}{#1}}
   \ITEE{#3}{WRudin1991}{
      \BIb{#2}{W. Rudin}
         {Functional Analysis}
         {McGraw-Hill Science}{1991}{#1}}
   \ITEE{#3}{TSaito1972}{
      \BiB{#2}{T. Sait\^{o}}{Generations of von Neumann algebras}
         {Lecture Notes in Math. vol. 247}{\textup{(}Lecture on Operator Algebras\textup{)}}
         {Springer, Berlin-Heidelberg-New York}{1972}{435--531}{#1}}
   \ITEE{#3}{KSakai,MYaguchi2003}{
      \BIB{#2}{K. Sakai and M. Yaguchi}
         {Characterizing manifolds modeled on certain dense subspaces of non-separable Hilbert spaces}
         {\jRN{TsukubaJM}}{27}{2003}{143--159}{#1}}
   \ITEE{#3}{SSakai1971}{
      \BIb{#2}{S. Sakai}
         {$\CCc^*$-Algebras and $\WWw^*$-Algebras}
         {Springer-Verlag, Berlin-Heidelberg-New York}{1971}{#1}}
   \ITEE{#3}{RSchori1971}{
      \BIB{#2}{R. Schori}
         {Topological stability for infinite-dimensional manifolds}
         {\jRN{ComposM}}{23}{1971}{87--100}{#1}}
   \ITEE{#3}{JTSchwartz1967}{
      \BIb{#2}{J.T. Schwartz}
         {$\WWw^*$-algebras}
         {Gordon and Breach, Science Publishers Inc., New York-London-Paris}{1967}{#1}}
   \ITEE{#3}{ZSemadeni1971}{
      \BIb{#2}{Z. Semadeni}
         {Banach Spaces of Continuous Functions (Vol. I)}
         {\jRN{PWN}}{1971}{#1}}
   \ITEE{#3}{JPSerre1951}{
      \BIB{#2}{J.-P. Serre}
         {Homologie singuli\`{e}re des espaces fibr\'{e}s}
         {\jRN{AnnM}}{54}{1951}{425--505}{#1}}
   \ITEE{#3}{DSherman2007}{
      \BIB{#2}{D. Sherman}
         {On the dimension theory of von Neumann algebras}
         {\jRN{MScand}}{101}{2007}{123--147}{#1}}
   \ITEE{#3}{WSierpinski1928}{
      \BIB{#2}{W. Sierpi\'{n}ski}
         {Sur les projections des ensembles compl\'{e}mentaires aux ensembles \textup{(A)}}
         {\jRN{FM}}{11}{1928}{117--122}{#1}}
   \ITEE{#3}{MSlocinski1980}{
      \BIB{#2}{M. S\l{}oci\'{n}ski}
         {On the Wold-type decomposition of a pair of commuting isometries}
         {\jRN{APM}}{37}{1980}{255--262}{#1}}
   \ITEE{#3}{RCSteinlage1975}{
      \BIB{#2}{R.C. Steinlage}
         {On Haar measure in locally compact $T_2$ spaces}
         {\jRN{AmJM}}{97}{1975}{291--307}{#1}}
   \ITEE{#3}{JStochel,FHSzafraniec1989}{
      \BIB{#2}{J. Stochel and F.H. Szafraniec}
         {On normal extensions of unbounded operators. III. Spectral properties}
         {\jRN{PublRIMSKyoto}}{25}{1989}{105--139}{#1}}
   \ITEE{#3}{JStochel,FHSzafraniec1989a}{
      \BIB{#2}{J. Stochel and F.H. Szafraniec}
         {The normal part of an unbounded operator}
         {\jRN{ProcKonink}}{92}{1989}{495--503}{#1}}
   \ITEE{#3}{AHStone1962}{
      \BIB{#2}{A.H. Stone}
         {Absolute $\FFf_{\sigma}$-spaces}
         {\jRN{PAMS}}{13}{1962}{495--499}{#1}}
   \ITEE{#3}{AHStone1962a}{
      \BIB{#2}{A.H. Stone}
         {Non-separable Borel sets}
         {\jRN{DissM}}{28}{1962}{41 pages}{#1}}
   \ITEE{#3}{AHStone1972}{
      \BIB{#2}{A.H. Stone}
         {Non-separable Borel sets II}
         {\jRN{GTopA}}{2}{1972}{249--270}{#1}}
   \ITEE{#3}{MHStone1937}{
      \BIB{#2}{M.H. Stone}
         {Application of the theory of Boolean rings to general topology}
         {\jRN{TAMS}}{41}{1937}{375--481}{#1}}
   \ITEE{#3}{MHStone1948}{
      \BIB{#2}{M.H. Stone}
         {The generalized Weierstrass approximation theorem}
         {\jRN{MMag}}{21}{1948}{167--184}{#1}}
   \ITEE{#3}{BSz-Nagy1947}{
      \BIB{#2}{B. Sz.-Nagy}
         {On uniformly bounded linear transformations in Hilbert space}
         {\jRN{ActaSM}}{11}{1947}{152--157}{#1}}
   \ITEE{#3}{WTakahashi1970}{
      \BIB{#2}{W. Takahashi}
         {A convexity in metric space and nonexpansive mappings, I}
         {\jRN{KodaiMSemRep}}{22}{1970}{142--149}{#1}}
   \ITEE{#3}{MTakesaki2002}{
      \BIb{#2}{M. Takesaki}
         {Theory of Operator Algebras I \textup{(Encyclopaedia of Mathematical Sciences, Volume 124)}}
         {Springer-Verlag, Berlin-Heidelberg-New York}{2002}{#1}}
   \ITEE{#3}{MTakesaki2003}{
      \BIb{#2}{M. Takesaki}
         {Theory of Operator Algebras II \textup{(Encyclopaedia of Mathematical Sciences, Volume 125)}}
         {Springer-Verlag, Berlin-Heidelberg-New York}{2003}{#1}}
   \ITEE{#3}{MTakesaki2003a}{
      \BIb{#2}{M. Takesaki}
         {Theory of Operator Algebras III \textup{(Encyclopaedia of Mathematical Sciences, Volume 127)}}
         {Springer-Verlag, Berlin-Heidelberg-New York}{2003}{#1}}
   \ITEE{#3}{TYTam1999}{
      \BIB{#2}{T.Y. Tam}
         {An extension of a result of Lewis}
         {\jRN{ELA}}{5}{1999}{1--10}{#1}}
   \ITEE{#3}{TYTam2000}{
      \BIB{#2}{T.Y. Tam}
         {Group majorization, Eaton triples and numerical range}
         {\jRN{LMLA}}{47}{2000}{11--28}{#1}}
   \ITEE{#3}{TYTam2002}{
      \BIB{#2}{T.Y. Tam}
         {Generalized Schur-concave functions and Eaton triples}
         {\jRN{LMLA}}{50}{2002}{113--120}{#1}}
   \ITEE{#3}{TYTam,WCHill2001}{
      \BIB{#2}{T.Y. Tam and W.C. Hill}
         {On $G$-invariant norms}
         {\jRN{LAA}}{331}{2001}{101--112}{#1}}
   \ITEE{#3}{AFTiman,IAVestfrid1983}{
      \BIB{#2}{A.F. Timan and I.A. Vestfrid}
         {Any separable ultrametric space can be isometrically imbedded in $l_2$}
         {\jRN{FAA}}{17}{1983}{70--71}{#1}}
   \ITEE{#3}{JTomiyama1958}{
      \BIB{#2}{J. Tomiyama}
         {Generalized dimension function for $\WWw^*$-algebras of infinite type}
         {\jRN{TohokuMJ} (2)}{10}{1958}{121--129}{#1}}
   \ITEE{#3}{HTorunczyk1970}{
      \BIB{#2}{H. Toru\'{n}czyk}
         {Remarks on Anderson's paper ``On topological infinite deficiency''}
         {\jRN{FM}}{66}{1970}{393--401}{#1}}
   \ITEE{#3}{HTorunczyk1970a}{
      \BIb{#2}{H. Toru\'{n}czyk}
         {$G$-$K$-absorbing and skeletonized sets in metric spaces}
         {Ph.D. thesis, Inst. Math. Polish Acad. Sci., Warszawa}{1970}{#1}}
   \ITEE{#3}{HTorunczyk1972}{
      \BIB{#2}{H. Toru\'{n}czyk}
         {A short proof of Hausdorff's theorem on extending metrics}
         {\jRN{FM}}{77}{1972}{191--193}{#1}}
   \ITEE{#3}{HTorunczyk1974}{
      \BIB{#2}{H. Toru\'{n}czyk}
         {Absolute retracts as factors of normed linear spaces}
         {\jRN{FM}}{86}{1974}{53--67}{#1}}
   \ITEE{#3}{HTorunczyk1975}{
      \BIB{#2}{H. Toru\'{n}czyk}
         {On Cartesian factors and the topological classification of linear metric spaces}
         {\jRN{FM}}{88}{1975}{71--86}{#1}}
   \ITEE{#3}{HTorunczyk1978}{
      \BIB{#2}{H. Toru\'{n}czyk}
         {Concerning locally homotopy negligible sets and characterization of $l_2$-manifolds}
         {\jRN{FM}}{101}{1978}{93--110}{#1}}
   \ITEE{#3}{HTorunczyk1980}{
      \BiB{#2}{H. Toru\'{n}czyk}{Characterization of infinite-dimensional manifolds}{in:}
         {Proceedings of the International Conference on Geometric Topology (Warsaw, 1978)}
         {\jRN{PWN}}{1980}{431--437}{#1}}
   \ITEE{#3}{HTorunczyk1981}{
      \BIB{#2}{H. Toru\'{n}czyk}
         {Characterizing Hilbert space topology}
         {\jRN{FM}}{111}{1981}{247--262}{#1}}
   \ITEE{#3}{HTorunczyk1985}{
      \BIB{#2}{H. Toru\'{n}czyk}
         {A correction of two papers concerning Hilbert manifolds}
         {\jRN{FM}}{125}{1985}{89--93}{#1}}
   \ITEE{#3}{KTsuda1985}{
      \BIB{#2}{K. Tsuda}
         {A note on closed embeddings of finite dimensional metric spaces}
         {\jRN{BLondMS}}{17}{1985}{273--278}{#1}}
   \ITEE{#3}{PSUrysohn1925}{
      \BIB{#2}{P.S. Urysohn}
         {Sur un espace m\'{e}trique universel}
         {\jRN{CRASParis}}{180}{1925}{803--806}{#1}}
   \ITEE{#3}{PSUrysohn1927}{
      \BIB{#2}{P.S. Urysohn}
         {Sur un espace m\'{e}trique universel}
         {\jRN{BullSM}}{51}{1927}{43--64, 74--96}{#1}}
   \ITEE{#3}{VVUspenskij1986}{
      \BIB{#2}{V.V. Uspenskij}
         {A universal topological group with a countable basis}
         {\jRN{FAA}}{20}{1986}{86--87}{#1}}
   \ITEE{#3}{VVUspenskij1990}{
      \BIB{#2}{V.V. Uspenskij}
         {On the group of isometries of the Urysohn universal metric space}
         {\jRN{CMUC}}{31}{1990}{181--182}{#1}}
   \ITEE{#3}{VVUspenskij2004}{
      \BIB{#2}{V.V. Uspenskij}
         {The Urysohn universal metric space is homeomorphic to a Hilbert space}
         {\jRN{TopA}}{139}{2004}{145--149}{#1}}
   \ITEE{#3}{VVUspenskij2008}{
      \BIB{#2}{V.V. Uspenskij}
         {On subgroups of minimal topological groups}
         {\jRN{TopA}}{155}{2008}{1580--1606}{#1}}
   \ITEE{#3}{VSVaradarajan1963}{
      \BIB{#2}{V.S. Varadarajan}
         {Groups of automorphisms of Borel spaces}
         {\jRN{TAMS}}{109}{1963}{191--220}{#1}}
   \ITEE{#3}{AMVershik1998}{
      \BIB{#2}{A.M. Vershik}
         {The universal Urysohn space, Gromov's metric triples, and random metrics on the series of natural numbers}
         {\jRN{UspekhiMN}}{53}{1998}{57--64}{#1} English translation: \jRN{RussMS}{} \textbf{53} (1998), 921--928. 
         Correction: \jRN{UspekhiMN}{} \textbf{56} (2001), p. 207. English translation: \jRN{RussMS}{} \textbf{56} 
         (2001), p. 1015.}
   \ITEE{#3}{AMVershik2002}{
      \BIb{#2}{A.M. Vershik}
         {Random metric spaces and the universal Urysohn space}
         {Fundamental Mathematics Today. 10th anniversary of the Independent Moscow University. MCCME Publ.}{2002}{#1}}
   \ITEE{#3}{NWeaver1999}{
      \BIb{#2}{N. Weaver}
         {Lipschitz Algebras}
         {World Scientific}{1999}{#1}}
   \ITEE{#3}{JWeidmann1980}{
      \BIb{#2}{J. Weidmann}
         {Linear Operators in Hilbert Spaces}
         {(Graduate Texts in Mathematics, vol. 68) Springer-Verlag New York Inc.}{1980}{#1}}
   \ITEE{#3}{JEWest1969}{
      \BIB{#2}{J.E. West}
         {Approximating homotopies by isotopies in Fr\'{e}chet manifolds}
         {\jRN{BAMS}}{75}{1969}{1254--1257}{#1}}
   \ITEE{#3}{JEWest1969a}{
      \BIB{#2}{J.E. West}
         {Fixed-point sets of transformation groups on infinite-product spaces}
         {\jRN{PAMS}}{21}{1969}{575--582}{#1}}
   \ITEE{#3}{JEWest1970}{
      \BIB{#2}{J.E. West}
         {The ambient homeomorphy of infinite-dimensional Hilbert spaces}
         {\jRN{PacJM}}{34}{1970}{257--267}{#1}}
   \ITEE{#3}{JHCWhitehead1949}{
      \BIB{#2}{J.H.C. Whitehead}
         {Combinatorial homotopy I}
         {\jRN{BAMS}}{55}{1949}{213--245}{#1}}
   \ITEE{#3}{GTWhyburn1942}{
      \BIb{#2}{G. T. Whyburn}
         {Analytic Topology}
         {Amer. Math. Soc. Colloquium Publications (vol. XXVIII), New York}{1942}{#1}}
   \ITEE{#3}{WWogen1969}{
      \BIB{#2}{W. Wogen}
         {On generators for von Neumann algebras}
         {\jRN{BAMS}}{75}{1969}{95--99}{#1}}
   \ITEE{#3}{RYTWong1967}{
      \BIB{#2}{R.Y.T. Wong}
         {On homeomorphisms of certain infinite dimensional spaces}
         {\jRN{TAMS}}{128}{1967}{148--154}{#1}}
   \ITEE{#3}{LYang,JZhang1987}{
      \BIB{#2}{L. Yang and J. Zhang}
         {Average distance constants of some compact convex space}
         {\jRN{JChinUST}}{17}{1987}{17--23}{#1}}
   \ITEE{#3}{PZakrzewski1993}{
      \BIB{#2}{P. Zakrzewski}
         {The existence of invariant $\sigma$-finite measures for a group of transformations}
         {\jRN{IsraelJM}}{83}{1993}{275--287}{#1}}
   \ITEE{#3}{PZakrzewski2002}{
      \BIb{#2}{P. Zakrzewski}
         {Measures on Algebraic-Topological Structures, Handbook of Measure Thoery}
         {E. Pap, ed., Elsevier, Amsterdam}{2002, 1091--1130}{#1}}
   \ITEE{#3}{KZhu2000}{
      \BIB{#2}{K. Zhu}
         {Operators in Cowen-Douglas classes}
         {\jRN{IllinoisJM}}{44}{2000}{767--783}{#1}}
   }

\newcommand{\mypaplist}[2][]{
   \ITEE{#2}{pn1}{
      \myBIB{Separate and joint similarity to families of normal operators}
         {\jRN[#1]{SM}}{149}{2002}{39--62}}
   \ITEE{#2}{pn2}{
      \myBIB{Locally arcwise connected metrizable spaces with the fixed point property are complete-metrizable}
         {\jRN[#1]{TopA}}{153}{2006}{1639--1642}}
   \ITEE{#2}{pn3}{
      \myBIB{Invariant measures for equicontinuous semigroups of continuous transformations of a compact Hausdorff space}
         {\jRN[#1]{TopA}}{153}{2006}{3373--3382}}
   \ITEE{#2}{pn4}{
      \myBIB{Approximation of the Hausdorff distance by the distance of continuous surjections}
         {\jRN[#1]{TopA}}{154}{2007}{655--664}}
   \ITEE{#2}{pn5}{
      \myBIB{Generalized Haar integral}
         {\jRN[#1]{TopA}}{155}{2008}{1323--1328}}
   \ITEE{#2}{pn6}{
      \myBIB{Integration and Lipschitz functions}
         {\jRN[#1]{RCMP}}{57}{2008}{391--399}}
   \ITEE{#2}{pn7}{
      \myBIB{Canonical Banach function spaces generated by Urysohn universal spaces. Measures as Lipschitz maps}
         {\jRN[#1]{SM}}{192}{2009}{97--110}}
   \ITEE{#2}{pn8}{
      \myBIB{Urysohn universal spaces as metric groups of exponent $2$}
         {\jRN[#1]{FM}}{204}{2009}{1--6}}
   \ITEE{#2}{pn9}{
      \myBIB{Central subsets of Urysohn universal spaces}
         {\jRN[#1]{CMUC}}{50}{2009}{445--461}}
   \ITEE{#2}{pn10}{
      \myBIB[P. Niemiec and T.Y. Tam]{A representation of $G$-in\-variant norms for Eaton triple}
         {\jRN[#1]{JCA}}{18}{2011}{59--65}}
   \ITEE{#2}{pn11}{
      \myBIB{Functor of extension of contractions on Urysohn universal spaces}
         {\jRN[#1]{ACS}}{}{2009}{\texttt{DOI: 10.1007/s10485-009-9218-z}}}
   \ITEE{#2}{pn12}{
      \myBIB{Ultra-$\mM$-separability}
         {\jRN[#1]{TopA}}{157}{2010}{669--673}}
   \ITEE{#2}{pn13}{
      \myBIB{Functor of extension of $\Lambda$-isometric maps between central subsets 
         of the unbounded Urysohn universal space}{\jRN[#1]{CMUC}}{51}{2010}{541--549}}
   \ITEE{#2}{pn14}{
      \myBIB{Normed topological pseudovector groups}{\jRN[#1]{ACS}}{}{2010}
         {\ITE{\equal{#1}{}}{\texttt{DOI: 10.1007/s10485\-010-9239-7}}{\texttt{DOI: 10.1007/s10485-010-9239-7}}}}
   \ITEE{#2}{pn15}{
      \myBIB{Topological structure of Urysohn universal spaces}
         {\jRN[#1]{TopA}}{158}{2011}{352--359}}
   \ITEE{#2}{pn16}{
      \myBIB{A note on invariant measures}
         {\jRN[#1]{OpusM}}{31}{2011}{425--431}}
   \ITEE{#2}{pn17}{
      \myBIB{Strengthened Stone-Weierstrass type theorem}
         {\jRN[#1]{OpusM}}{31}{2011}{645--650}}
   \ITEE{#2}{pnX2}{
      \myBAPP{Functor of continuation in Hilbert cube and Hilbert space}
         {to appear in \jRN[#1]{FM}}}
   \ITEE{#2}{pnX3}{
      \myBAPP{Norm closures of orbits of bounded operators}
         {to appear.}}
   \ITEE{#2}{pnX6}{
      \myBAPP{Extending maps by injective $\sigma$-$Z$-maps in Hilbert manifolds}
         {to appear in \jRN[#1]{BullPol}}}
   \ITEE{#2}{pnX7}{
      \myBAPP{Spaces of measurable functions}
         {submitted to \jRN[#1]{CollectM}}}
   \ITEE{#2}{pnX8}{
      \myBAPP{Normal systems over ANR's, rigid embeddings and nonseparable absorbing sets}
         {submitted to \jRN[#1]{ActaMSinES}}}
   \ITEE{#2}{pnX9}{
      \myBAPP{Borel structure of the spectrum of a closed operator}
         {submitted to \jRN[#1]{SM}}}
   \ITEE{#2}{pnX10}{
      \myBAPP{Central points and measures and dense subsets of compact metric spaces}
         {submitted to \jRN[#1]{TopMethNA}}}
   \ITEE{#2}{pnX11}{
      \myBAPP{Generalized absolute values and polar decompositions of a bounded operator}
         {submitted to \jRN[#1]{IEOT}.}}
   \ITEE{#2}{pnX12}{
      \myBAPP{Ultrametrics, extending of Lipschitz maps and nonexpansive selections}
         {accepted for publication in \jRN[#1]{HJM}}}
   \ITEE{#2}{pnX13}{
      \myBAPP{A note on ANR's}
         {submitted to \jRN[#1]{TopA}}}
   \ITEE{#2}{pnX14}{
      \myBAPP{Problem with almost everywhere equality}
         {submitted to \jRN[#1]{ArchM}}}
   \ITEE{#2}{pnX15}{
      \myBAPP{Universal valued Abelian groups}
         {submitted to \jRN[#1]{LNM}}}
   \ITEE{#2}{pnX16}{
      \myBAPP{Unitary equivalence and decompositions of finite systems of closed densely defined operators 
         in Hilbert spaces}{submitted to \jRN[#1]{DissM}}}
   }

\begin{document}

\title[Extending maps in Hilbert manifolds]{Extending maps by injective $\sigma$-$Z$-maps\\in Hilbert manifolds}
\myData
\begin{abstract}
The aim of the paper is to prove that if $M$ is a metrizable manifold modelled on a Hilbert space of dimension $\alpha
\geqsl \aleph_0$ and $F$ is its $\sigma$-$Z$-set, then for every completely metrizable space $X$ of weight no greater
than $\alpha$ and its closed subset $A$, for any map $f \dd X \to M$, each open cover $\UUu$ of $M$ and a sequnce
$(A_n)_n$ of closed subsets of $X$ disjoint from $A$ there is a map $g\dd X \to M$ $\UUu$-homotopic to $f$ such that
$g\bigr|_A = f\bigr|_A$, $g\bigr|_{A_n}$ is a closed embedding for each $n$ and $g(X \setminus A)$ is a $\sigma$-$Z$-set
in $M$ disjoint from $F$. It is shown that if $f(\partial A)$ is contained in a locally closed $\sigma$-$Z$-set in $M$
or $f(X \setminus A) \cap \overline{f(\partial A)} = \varempty$, the map $g$ may be taken so that
$g\bigr|_{X \setminus A}$ be an embedding. If, in addition, $X \setminus A$ is a connected manifold modelled on the same
Hilbert space as $M$ and $\overline{f(\partial A)}$ is a $Z$-set in $M$, then there is a $\UUu$-homotopic to $f$ map
$h\dd X \to M$ such that $h\bigr|_A = f\bigr|_A$ and $h\bigr|_{X \setminus A}$ is an open embedding.
\end{abstract}
\maketitle


In \cite{west}, West has proved that any homotopy $(f_t\dd X \to M)_{t \in I}$ from a separable completely metrizable
space $X$ into a separable metrizable manifold $M$ modelled on an infinite-dimensional Fr\'{e}chet space can be
approximated by homotopies $(h_t\dd X \to M)_{t \in I}$ such that $h_t = f_t$ for $t=0,1$ and each $h_t$
with $t \in (0,1)$ is a closed embedding. He has also made a note that the homotopy $(h_t)_{t \in I}$ may be modified
in such a way that it be an embedding of $X \times (0,1)$ into $M$. This claim was applied by Anderson and McCharen
in their joint paper \cite{and-mcch} on extending homeomorphisms between $Z$-set of (separable) Fr\'{e}chet manifolds.
Unfortunately, this statement of West is incorrect (see \EXM{false} below). In this paper we try to give sufficient
conditions under which a map of a complete metric space $X$ into a Hilbert manifold $M$ could be approximated,
in the limitation topology and relative a given closed set $A \subset X$, by mappings $g$ which are embeddings
on $X \setminus A$. Moreover, we request that $g$ is a $Z$-embedding on any of the countably many closed subsets
of $X \setminus A$, or, when $X \setminus A$ is a manifold modelled on the same Hilbert space as $M$,
that $g\bigr|_{X \setminus A}$ be an open embedding of $X \setminus A$ into $M$. Our proofs totally differ
from that of West \cite{west} and depend on the argument used by Toru\'{n}czyk in \cite[Proof of~3.1]{tor2}.\par
The paper is organized as follows. In Section~1 we establish notation and terminology and cite theorems which shall be
applied in the next sections. The second part is devoted to the proof of the main lemma on extending maps which take
values in ANR's. Section~3 deals with extending maps from completely metrizable spaces into infinite-dimensional
Hilbert manifolds by injections whose images are $\sigma$-$Z$-sets. In the last part we give conditions under which
a map is extendable by an embedding.

\SECT{Theorems to quote}

In this paper $I$ denotes the unit interval $[0,1]$. The letters $X$, $Y$, $Z$, $K$, etc. stand for metrizable spaces.
By a \textit{map} we mean a continuous function. Whenever $g$ is a function, $\im g$ stands for the image of $g$.
If, in addition, $g$ takes values in a topological space, by $\overline{\im}\,g$ we denote the closure of the image
in the whole space. If $A$ is a subset of $X$, $\intt A$, $\bar{A}$ and $\partial A$ stand for, respectively,
the interior, the closure and the boundary of $A$ in the whole space $X$. If $B \subset X$ is a superset of $A$,
by $\intt_B A$, $\cll_B A$ and $\partial_B A$ we denote the interior, closure and boundary of $A$ relative to $B$.
We use $w(X)$ to denote the topological weight of $X$.\par
Whenever $Y$ is a metrizable space, $\cov(Y)$ and $\Metr(Y)$ denote the collections of all open covers of $Y$ and
of all metrics on $Y$ inducing its topology, respectively. For every $\varrho \in \Metr(Y)$, $B_{\varrho}(y,r)$
and $\bar{B}_{\varrho}(y,r)$ stand for the open and the closed $\varrho$-ball (respectively) in $Y$ with center
at $y \in Y$ and of radius $r > 0$. Following Toru\'{n}czyk \cite{tor1}, for a map $f$ of a metrizable space $X$ into $Y$,
$B(f,\UUu)$ with $\UUu \in \cov(Y)$ consists of all maps $g\dd X \to Y$ which are \textit{$\UUu$-close} to $f$, that is,
$g$ belongs to $B(f,\UUu)$ iff for every $x \in X$ there is $U \in \UUu$ such that $\{f(x),g(x)\} \subset U$. Similarly,
for $\varrho \in \Metr(Y)$ and a map $\alpha\dd Y \to (0,+\infty)$, $B_{\varrho}(f,\alpha)$ is the set of all maps
$g\dd X \to Y$ such that $\varrho(f(x),g(x)) \leqsl \alpha(f(x))$ for any $x \in X$. On the space $\CCc(X,Y)$ of all maps
of $X$ into $Y$ we always consider the \textit{limitation topology} in which $\{B(f,\UUu)\dd\ \UUu \in \cov(Y)\}$
(respectively $\{B_{\varrho}(f,\alpha)\dd\ \alpha \in \CCc(Y,(0,+\infty))\}$ with fixed $\varrho \in \Metr(Y)$) may serve
as a base of open (closed) neighbourhoods of a given map $f$. For basic properties of this topology the reader is referred
to \cite{tor1}, \cite{bowers}. The set of all [closed; open] embeddings of $X$ into $Y$ is denoted by $\Emb(X,Y)$
[$\Emb^c(X,Y)$; $\Emb^o(X,Y)$].\par
By a \textit{Hilbert manifold} we mean a metrizable space which admits an open cover by sets homeomorphic
to some infinite-dimensional Hilbert space. Any metrizable manifold modelled on an infinite-dimensional Fr\'{e}chet
space is a Hilbert manifold (since every Fr\'{e}chet space is homeomorphic to a Hilbert one --- see \cite{tor1,tor2}).
Hilbert manifolds are completely metrizable ANR's. For simplicity, we say that a space is an \textit{$\alpha$-manifold},
where $\alpha$ is an infinite cardinal, if it is a manifold modelled on the Hilbert space of dimension $\alpha$.\par
If $\UUu$ is any collection of subsets of $Y$ ($\UUu$ need not cover $Y$ and the members of $\UUu$ need not be
open) and $F\dd X \times I \to Y$ is a homotopy, then $F$ is said to be a \textit{$\UUu$-homotopy} iff for any $x \in X$
either $F(\{x\} \times I)$ consists of a single point or there is $U \in \UUu$ such that $F(\{x\} \times I) \subset U$.
Two maps $f,g\dd X \to Y$ are \textit{$\UUu$-homotopic in $Y$} if there is a $\UUu$-homotopy $X \times I \to Y$
which connects $f$ and $g$. Additionally, if $f\bigr|_A = g\bigr|_A$, $f$ and $g$ are said to be
\textit{$\UUu$-homotopic in $Y$ relative $A$} if there is a $\UUu$-homotopy $F$ connecting $f$ and $g$ such that
$F(a,t) = f(a)$ for every $a \in A$ and $t \in [0,1]$.\par
The following is known as Toru\'{n}czyk's Lemma (see \cite[Lemma~1.1]{tor1}, and \cite{bowers} for proof).

\begin{lem}{tor}
Let $Y$ be completely metrizable, $F$ a subspace of $\CCc(X,Y)$ and $U_1,U_2,\ldots$ open subsets of $\CCc(X,Y)$.
If $U_n \cap F$ is dense in $F$ for each $n \geqsl 1$, then maps in $F$ are approximable by elements of $F_{\varrho}
\cap \bigcap_{n=1}^{\infty} U_n$, where $F_{\varrho}$ denotes the closure of $F$ in the topology of $\varrho$-uniform
convergence in $\CCc(X,Y)$ and $\varrho \in \Metr(Y)$.
\end{lem}

In the sequel we shall also apply the next two results of Toru\'{n}czyk and the one of Henderson and Schori.

\begin{lemm}{ANR-open}{\mbox{\cite[Lemma~1.3]{tor1}}}
If $Y$ is an ANR and $A$ is a closed subset of $X$, then the map $\CCc(X,Y) \ni f \mapsto f\bigr|_A \in \CCc(A,Y)$
is open.
\end{lemm}

\begin{lemm}{emb-dense}{\mbox{\cite{tor1}}}
If $X$ is completely metrizable and of weight no greater than $\alpha \geqsl \aleph_0$ and $M$ is an $\alpha$-manifold,
then the set $\Emb^c(X,M)$ is dense in $\CCc(X,Y)$.
\end{lemm}

\begin{thmm}{open-emb-dense}{\mbox{\cite{h-s}}}
If $M$ and $N$ are two $\alpha$-manifolds (where $\alpha$ is an infinite cardinal) and $M$ has no more than
$\alpha$ components (i.e. $w(M) = \alpha$), then the set $\Emb^o(M,N)$ is dense in $\CCc(M,N)$.
\end{thmm}

The next two lemmas are easy to prove.

\begin{lem}{comp-proj}
Let $K$ be compact and $p\dd Y \times K \to Y$ be the natural projection. Then the map $\CCc(X,Y \times K) \ni f
\mapsto p \circ f \in \CCc(X,Y)$ is open.
\end{lem}

\begin{lem}{closed-agree}
If $F$ is a closed subset of $Y$, then the limitation topology of $\CCc(X,F)$ coincides with the one induced
by the limitation topology of $\CCc(X,Y)$, when $\CCc(X,F)$ is considered as a subset of $\CCc(X,Y)$.
\end{lem}

Recall that a subset $P$ of $Y$ is \textit{locally closed} (in $Y$) if every point $p$ of $P$ has an open in $Y$
neighbourhood $U$ such that $P \cap U$ is relatively closed in $U$. Equivalently, $P$ is locally closed iff
$\bar{P} \setminus P$ is closed in $Y$.

\begin{lem}{hom-ext}
Let $A$ be a closed subset of $X$ and $P$ a locally closed ANR-set in $Y$, and let a map $f\dd X \to Y$ be such that
$f(X \setminus A) \subset P$. Then, for every $\UUu \in \cov(P)$ there exists a neighbourhood $G$
of $f\bigr|_{X \setminus A}$ in $\CCc(X \setminus A,P)$ such that each map $g\dd X \to Y$ satisfying
$g\bigr|_A = f\bigr|_A$ and $g\bigr|_{X \setminus A} \in G$ is $\UUu$-homotopic to $f$ in $Y$ relative $A$.
\end{lem}
\begin{proof}
It is well known that ANR's are the so-called \textit{locally equiconnected spaces} (\cite{dug}, \cite{fox}, \cite{serre}),
that is, there is an open subset $\Omega$ of $P \times P$ containing the diagonal and a map $\lambda: \Omega \times [0,1]
\to P$ such that $\lambda(x,y,0) = x$, $\lambda(x,y,1) = y$ and $\lambda(z,z,t) = x$ for every $(x,y) \in \Omega$,
$z \in P$ and $t \in [0,1]$. Let $d \in \Metr(Y)$. If $P$ is closed in $Y$, put $\theta \equiv 1$
(as a member of $\CCc(P,(0,\infty))$). Otherwise let $\theta: P \ni z \mapsto \frac12 \dist_d(z,\bar{P} \setminus P)
\in (0,+\infty)$. Now every $z \in P$ has an open in $P$ neighbourhood $V_z$ such that
\begin{equation}\label{eqn:aux}
V_z \times V_z \subset \Omega \textup{ and } \lambda(V_z \times V_z \times [0,1]) \subset U \cap B_d(z,\theta(z))
\end{equation}
for some $U \in \UUu$. Put $\VVv = \{V_z\dd\ z \in P\} \in \cov(P)$ and $G = B(f\bigr|_{X \setminus A},\VVv) \subset
\CCc(X \setminus A, P)$ and assume that $g\dd X \to Y$ is as in the statement of the lemma. Define $F\dd X \times [0,1]
\to Y$ by $F(x,t) = \lambda(f(x),g(x),t)$ for $x \in X \setminus A$ and $F(x,t) = f(x)$ for $x \in A$. By \eqref{eqn:aux},
$F$ is well defined and we only need to check that it is continuous at points of $A$. Fix a sequence $(x_n,t_n) \in
(X \setminus A) \times [0,1]$ convergent to $(a,t) \in A \times [0,1]$. If $f(a) \in P$, then $(f(x_n),g(x_n),t_n) \to
(f(a),f(a),t)\ (n \to \infty)$ and the continuity of $\lambda$ gives $\lim_{n\to\infty} F(x_n,t_n) = f(a)$. Now suppose
that $f(a) \notin P$. Since $f(x_n) \in P$, we see that $f(a) \in \bar{P} \setminus P$. Take $z_n \in P$ for which $g(x_n),
f(x_n) \in V_{z_n}$. Then, by \eqref{eqn:aux}, both $F(x_n,t_n)$ and $f(x_n)$ belong to $B_d(z_n,\theta(z_n))$,
so $d(F(x_n,t_n),f(x_n)) \leqsl 2 \theta(z_n)$. But $\theta(z_n) \leqsl \frac12 d(z_n,f(a)) \leqsl \frac12 d(z_n,f(x_n)) +
\frac12 d(f(x_n),f(a)) \leqsl \frac12 \theta(z_n) + \frac12 d(f(x_n),f(a))$ and thus $\theta(z_n) \leqsl d(f(x_n),f(a))$.
Finally we obtain $d(F(x_n,t_n),f(a)) \leqsl d(F(x_n,t_n),f(x_n)) + d(f(x_n),f(a)) \leqsl 3 d(f(x_n),f(a)) \to 0$.
\end{proof}

\SECT{Main lemma}

The proof of the following result is based on the Proof of~3.1 given by Toru\'{n}czyk in \cite{tor2}.

\begin{lem}{main}
Let $P$ be a completely metrizable ANR-set in a metric space $(Y,\varrho)$ and $A, A_1, A_2, A_3, \ldots$ be closed
subsets of $X$ such that $A \cap A_n = \varempty \neq A_n$ for any $n \in \NNN$. Suppose also that, $G_n$ is a dense
$\GGg_{\delta}$-set in $\CCc(A_n,P)$ ($n \in \NNN$) and that a map $\lambda\dd X \setminus A \to (0,\infty)$ satisfies
$\inf \lambda(A_n) > 0$ for all $n \in \NNN$. Then, given a mapping $f\dd X \to Y$ with $f(X \setminus A) \subset P$,
and a neighbourhood $G$ of $f\bigr|_{X \setminus A}$ in $\CCc(X \setminus A,P)$, there is a mapping $h\dd X \to Y$
satisfying the following conditions:
\begin{enumerate}[\upshape({A}1)]
\item $h\bigr|_A = f\bigr|_A$, $h(X \setminus A) \subset P$,
\item $h\bigr|_{A_n} \in \GGg_n$ for all $n \in \NNN$,
\item $h\bigr|_{X \setminus A} \in G$ and $\varrho(h(x),f(x)) < \lambda(x)$ for all $x \in X \setminus A$.
\end{enumerate}
\end{lem}
\begin{proof}
Take maps $\delta\dd X \to I$ and $\theta\dd P \to (0,1]$ such that
\begin{equation}\label{eqn:aux100}
B_{\varrho}(f\bigr|_{X \setminus A},\theta) \subset G
\end{equation}
and
\begin{equation}\label{eqn:delta}
\delta \leqsl \lambda, \qquad \delta^{-1}(\{0\}) = A, \qquad \inf \delta(A_n) > 0 \textup{ for each } n \in \NNN.
\end{equation}
Let $d$ be the maximum metric on the space $P \times (0,1]$, that is, $$d((p,s),(q,t)) = \max(\varrho(p,q),|s-t|).$$
There is an open subset $U$ of $\CCc(X \setminus A,P \times (0,1])$ such that
\begin{multline}\label{eqn:aux3}
(f\bigr|_{X \setminus A},\delta\bigr|_{X \setminus A}) \in U \subset \{F \in \CCc(X \setminus A,P \times (0,1])\dd\\
d(F(x),(f(x),\delta(x))) \leqsl \frac12 \min(\theta(f(x)),\delta(x))\ (x \in X \setminus A)\}.
\end{multline}
Fix $n \geqsl 1$. By \LEM{ANR-open}, the set $U\bigr|_{A_n} := \{F\bigr|_{A_n}\dd\ F \in U\}$ is open
in $\CCc(A_n,P \times (0,1])$ and the map
$$
\varphi_n\dd U \ni F \mapsto F\bigr|_{A_n} \in U\bigr|_{A_n}
$$
is open as well. Let $m_n = \frac12 \inf \delta(A_n) > 0$. By \eqref{eqn:aux3}, $F(A_n) \subset P \times [m_n,1]$ for each
$F \in U$. This yields that $U\bigr|_{A_n} \subset \CCc(A_n,P \times [m_n,1])$. By \LEM{closed-agree}, $U\bigr|_{A_n}$
is open in $\CCc(A_n,P \times [m_n,1])$ and $\varphi_n$, as a map of $U$ into $\CCc(A_n,P \times [m_n,1])$, is open
as well. Now let $p\dd P \times (0,1] \to P$ be the natural projection and put $$\psi_n\dd U\bigr|_{A_n} \ni v \mapsto
p \circ v \in \CCc(A_n,P).$$ By \LEM{comp-proj}, $\psi_n$ is open. We conclude from this that the set
$D_n = \varphi_n^{-1}(\psi_n^{-1}(\GGg_n))$ is a dense $\GGg_{\delta}$ subset of $U$.\par
Now \LEM{tor} shows that the set $D = \bigcap_{n=1}^{\infty} D_n$ is dense in $U$. Take $u \in D \subset
\CCc(X \setminus A,P \times (0,1])$
and define $h\dd X \to Y$ as follows: $h\bigr|_A = f\bigr|_A$ and $h\bigr|_{X \setminus A} = p \circ u$. Thanks to
\eqref{eqn:aux3} and \eqref{eqn:delta}, $h$ is continuous. What is more, by construction, $h$ satisfies (A1) and (A2).
Finally, $h\bigr|_{X \setminus A} \in G$ because of \eqref{eqn:aux100} and \eqref{eqn:aux3} and the remainder of (A3)
follows from \eqref{eqn:aux3} and \eqref{eqn:delta} (note that $u \in U$).
\end{proof}

In the sequel we shall also need the next result. Since its proof is similar to that of \LEM{main} (but simpler),
we omit it.

\begin{lem}{main2}
Let $A$ and $P$ be subsets of $X$ and $Y$, respectively, such that $A$ is closed in $X$ and $P$ is an ANR.
Let $\SSs$ be a dense subset of $\CCc(X \setminus A,P)$. If $f \in \CCc(X,Y)$ is such that $f(X \setminus A) \subset P$
and $\overline{f(\partial A)} \cap P = \varempty$, then for every neighbourhood $G$ of $f\bigr|_{X \setminus A}$
in $\CCc(X \setminus A, P)$ there is a map $h\dd X \to Y$ such that $h\bigr|_A = f\bigr|_A$ and $h\bigr|_{X \setminus A}
\in G \cap \SSs$.
\end{lem}

We end the section with the following

\begin{exm}{counter}
As this example shows, the assumption of \LEM{main} that $\inf \lambda(A_n) > 0$ for each $n$ cannot be
omitted. Let $X = [0,+\infty)$; $Y = P = I$; $A = \{0\}$; $A_1 = [1,+\infty)$; $\GGg_1 = \{g \in \CCc(A_1,P)\dd\
g(x) \not\to 0\ (x \to \infty)\}$; $f\dd X \to Y$, $f \equiv 0$; $\lambda\dd X \to I$, $\lambda(t) = 1$
for $t \leqsl 1$ and $\lambda(t) = \frac1t$ for $t \geqsl 1$; and $\UUu = \{P\}$. Note that $\GGg_1$ is open and dense
in $\CCc(X,Y)$ and there is no map $h\dd X \to Y$ such that $|h(x) - f(x)| \leqsl \lambda(x)$ for each $x \in X$
and $h\bigr|_{A_1} \in \GGg_1$.
\end{exm}

\SECT{Extending maps by injective $\sigma$-$Z$-maps}

We begin with

\begin{dfn}{small}
For a subset $B$ of $Y$, let $\SsS_Y(X,B)$ be the collection of all maps $g\dd X \to Y$ such that $\overline{\im}\,g
\subset B$ (the closure taken in $Y$). Note that if $B$ is open or $\GGg_{\delta}$, so is $\SsS_Y(X,B)$.
\end{dfn}

Following Toru\'{n}czyk \cite{tor,tor0,tor1}, we say that a closed subset $A$ of $X$ is a \textit{$Z$-set} in $X$,
if the space $\CCc(Q,X \setminus A)$, where $Q$ denotes the Hilbert cube, is dense in $\CCc(Q,X)$.
(If $X$ is an ANR, this definition is equivalent to the original one by Anderson \cite{anderson}.) Similarly, $A$ is
said to be a \textit{strong $Z$-set} in $X$ iff for every $\UUu \in \cov(X)$ there is a map $u\dd X \to X$ which is
$\UUu$-close to the identity map $\id_X$ on $X$ and $A \cap \overline{\im}\,u = \varempty$
(cf. e.g. \cite{bbmw}, \cite{b-m}, \cite{dij1,dij2}). In other words, $A = \bar{A}$ is a $Z$-set in $X$
iff $\SsS_X(Q,X \setminus A)$ is dense in $\CCc(Q,X)$, and $A$ is a strong $Z$-set in $X$ iff $\SsS_X(X,X \setminus A)$
is dense in $\CCc(X,X)$. Countable unions of [strong] $Z$-sets are called [\textit{strong}] \textit{$\sigma$-$Z$-sets}.
If $X$ is complete metrizable and $B$ is its $\sigma$-$Z$-set, then the set $\SsS_X(Q,X \setminus B)$ is a dense
$\GGg_{\delta}$-set, and if $B$ is a strong $\sigma$-$Z$-set, the same is true with $Q$ replaced by $X$ (this follows
from Toru\'{n}czyk's Lemma).\par
Not every $Z$-set in an ANR is a strong $Z$-set (\cite[Key example, p.~56]{bbmw}). However, by a theorem
due to Henderson \cite{henderson}, $Z$-sets in Hilbert manifolds are strong $Z$-sets (for other results in this matter
see \cite{b-m}). This fact will be used by us several times.\par
We say that a map $f\dd X \to Y$ is a \textit{$Z$-map} [a \textit{$\sigma$-$Z$-map}] iff $\im f$ is a $Z$-set
[$\sigma$-$Z$-set] in $Y$. We similarly define $Z$-embeddings and $\sigma$-$Z$-embeddings (cf. \cite{tor1})
Note that $Z$-maps are closed. It is well known (\cite{tor1}) that if $X$ is completely metrizable and of weight
no greater than $\alpha \geqsl \aleph_0$ and $M$ is an $\alpha$-manifold, then $Z$-embeddings of $X$ into $M$
are dense in $\CCc(X,M)$. In the following result we only need to know that $Z$-maps are dense.

\begin{lem}{big-Z}
Every Hilbert manifold $M$ contains a $\sigma$-$Z$-set $F$ such that each closed subset of $M$ disjoint from $F$
is a $Z$-set in $M$.
\end{lem}
\begin{proof}
Take a sequence of $Z$-maps $f_n\dd M \to M$ which converges uniformly to $\id_M$ with respect to a fixed metric of $M$
and put $F = \bigcup_{n=1}^{\infty} \im f_n$.
\end{proof}

\begin{cor}{dense-Gd}
Let $M$ be an $\alpha$-manifold ($\alpha \geqsl \aleph_0$), $K$ a $\sigma$-$Z$-set in $M$ and let $X$ be a completely
metrizable space of weight no greater than $\alpha$. Then there is a dense $\GGg_{\delta}$-subset $G$ of $\CCc(X,M)$
which consists of $Z$-embeddings whose images are disjoint from $K$.
\end{cor}
\begin{proof}
Take $F$ as in \LEM{big-Z} and put $$G = \Emb^c(X,M) \cap \SsS_M(X,M \setminus (F \cup K)).$$ By Henderson's theorem
\cite{henderson}, $F \cup K$ is a strong $\sigma$-$Z$-set in $M$ and thus $\SsS_M(X,M \setminus (F \cup K))$ is a dense
$\GGg_{\delta}$-set in $\CCc(X,M)$. Thanks to \LEM{emb-dense}, also $\Emb^c(X,M)$ is a dense $\GGg_{\delta}$-set
in $\CCc(X,M)$. Finally, \LEM{tor} yields that $G$ is dense as well. The remainder of the assertion is clear.
\end{proof}

Now we are able to state and prove the main result of this section.

\begin{thm}{Z-maps}
Let $A$ be a closed subset of $X$ such that $X \setminus A$ is completely metrizable and $w(X \setminus A) \leqsl \alpha$
(where $\alpha \geqsl \aleph_0$) and let $M \subset Y$ be an $\alpha$-manifold. Let $A_1,A_2,\ldots$ be nonempty closed
subsets of $X$ disjoint from $A$ and let $B \subset M$ be a $\sigma$-$Z$-set in $M$. Let a map $\lambda\dd X \setminus A
\to (0,+\infty)$ be such that $\inf \lambda(A_n) > 0$ for each $n \in \NNN$. If $f\dd X \to Y$ is such a map that
$f(X \setminus A) \subset M$, then for every neighbourhood $G$ of $f\bigr|_{X \setminus A}$ in $\CCc(X \setminus A,M)$
and $\varrho \in \Metr(Y)$ there is a map $h\dd X \to Y$ satisfying the following conditions:
\begin{enumerate}[\upshape(Z1)]
\item $h\bigr|_A = f\bigr|_A$, $h(X \setminus A) \subset M \setminus B$,
\item $h\bigr|_{A_n}$ is a $Z$-embedding into $M$ for all $n \in \NNN$,
\item $h\bigr|_{X \setminus A} \in G$ and $\varrho(h(x),f(x)) < \lambda(x)$ for all $x \in X \setminus A$.
\end{enumerate}
\end{thm}
\begin{proof}
Enlarging the sets $A_n$, we may and do assume that
\begin{equation}\label{eqn:aux5}
X \setminus A = \bigcup_{n=1}^{\infty} A_n.
\end{equation}
Since $w(X \setminus A) \leqsl \alpha$ and $X \setminus A$ is completely metrizable, by \COR{dense-Gd}, for each
$n \in \NNN$ there is a dense $\GGg_{\delta}$-subset $\GGg_n$ of $\CCc(A_n,M)$ consisting of $Z$-embeddings whose images
are disjoint from $B$. Now putting $P = M$ and applying \LEM{main}, we obtain a map $h\dd X \to Y$ such that the conditions
(A1)--(A3) are fulfilled. This yields (Z3) and (Z2). Finally, (Z1) is also satisfied because of (A2) and \eqref{eqn:aux5}.
\end{proof}

\begin{rem}{stronger}
As the above proof shows, the map $h$ appearing in the statement of \THM{Z-maps} may be chosen in such a way that
it additionally fulfills the following condition:
\begin{enumerate}[\upshape(Z1)]\addtocounter{enumi}{3}
\item $h\bigr|_{X \setminus A}$ is a one-to-one $\sigma$-$Z$-map as a map of $X \setminus A$ into $M$ (this is guaranteed
   by \eqref{eqn:aux5} after enlarging the sets $A_n$).
\end{enumerate}
What is more, with use of \LEM{hom-ext}, $h$ may be forced to satisfy also:
\begin{enumerate}[\upshape(H1)]
\item the maps $h\bigr|_{X \setminus A}$ and $f\bigr|_{X \setminus A}$ are $\UUu$-homotopic in $M$ and
\item if $M$ is locally closed in $Y$, then $h$ and $f$ are $\UUu$-homotopic in $Y$
\end{enumerate}
where $\UUu$ is an arbitrarily given relatively open cover of $M$. Both the points (H1) and (H2)
may be added also in the statements of \PRO{emb} and \THM{open}. This observation shall be used in the sequel
(see Corollaries \ref{cor:homotopy}, \ref{cor:main}, \ref{cor:homo} and \ref{cor:open}).
\end{rem}

As a consequence of \THM{Z-maps} we get a generalization of West's theorem \cite{west} and a strengthened version
of \LEM{emb-dense}:

\begin{cor}{homotopy}
Let $X$ be completely metrizable with $w(X) \leqsl \alpha$, $M$ an $\alpha$-manifold and $B$ its $\sigma$-$Z$-set.
\begin{enumerate}[\upshape(a)]
\item For every homotopy $F\dd X \times I \to M$ and each $\UUu \in \cov(M)$ there is a homotopy $H\dd X \times I\to M$
   $\UUu$-close to $F$ such that $F(\cdot,t) = H(\cdot,t)$ for $t = 0,1$ and $H\bigr|_{X \times [\epsi,1-\epsi]}$
   is a $Z$-embedding with image disjoint from $B$ for any $\epsi \in (0,\frac12)$.
\item For every map $f\dd X \to M$ and each $\UUu \in \cov(M)$ there is a $\UUu$-homotopy $F\dd X \times I \to M$
   such that $F(\cdot,0) = f$ and for any $\epsi \in (0,1)$, $F\bigr|_{X \times [\epsi,1]}$ is a $Z$-embedding
   with image disjoint from $B$. In particular, $f_t = F(\cdot,t)$ (with $t \in (0,1]$) is a $Z$-embedding of $X$
   into $M$.
\end{enumerate}
\end{cor}

\SECT{Extending maps by embeddings}

In this section we give sufficient conditions under which an arbitrary map $f\dd X \to M$ is approximable by maps
$g\dd X \to M$ such that $g\bigr|_A = f\bigr|_A$ and $g\bigr|_{X \setminus A}$ is an embedding.\par
The proof of the following is left as a simple exercise.

\begin{lem}{disjoint}
If $h\dd X \to Y$ is such a map that $h\bigr|_U$ is an embedding, where $U$ is open in $X$,
then $\overline{h(\partial U)} \cap h(U) = \varempty$.
\end{lem}

The property stated in the above lemma forces us to make some restrictions on a map which we want to extend
by an embedding.\par
For need of the next result, note that a subset $F$ of a Hilbert manifold $M$ is a locally closed $\sigma$-$Z$-set
in $M$ iff $F$ is a $Z$-set in some open in $M$ neighbourhood of $F$.

\begin{pro}{emb}
Let $A$ be a closed subset of $X$ such that $X \setminus A$ is completely metrizable and of weight no greater
than $\alpha$ ($\alpha \geqsl \aleph_0$). Let $M \subset Y$ be an $\alpha$-manifold and $B$ its $\sigma$-$Z$-set.
If $f\dd X \to Y$ is such a map that there is a closed subset $K$ of $Y$ such that $f(X \setminus A) \subset
M \setminus K$ and $M \cap (\overline{f(\partial A)} \setminus K)$ is a $\sigma$-$Z$-set in $M$, then for each
neighbourhood $G$ of $f\bigr|_{X \setminus A}$ in $\CCc(X \setminus A,M)$ there is a map $h\dd X \to Y$ satisfying
the following conditions:
\begin{enumerate}[\upshape(E1)]
\item $h\bigr|_A = f\bigr|_A$, $h(X \setminus A) \subset M \setminus (B \cup K)$ and $h\bigr|_{X \setminus A} \in G$,
\item $h\bigr|_{X \setminus A}$ is a $\sigma$-$Z$-embedding of $X \setminus A$ into $M$ whose image is locally closed
   in $M$,
\item if $f(\partial A) \subset M \setminus K$ and $f\bigr|_{\partial A}$ is an embedding, then so is
   $h\bigr|_{\overline{X \setminus A}}$.
\end{enumerate}
\end{pro}
\begin{proof}
Put $M' = M \setminus K$ and
\begin{equation}\label{eqn:aux6}
B' = \overline{f(\partial A)} \cap M'.
\end{equation}
Then $M'$ is an $\alpha$-manifold, $B'$ its closed $\sigma$-$Z$-set and thus it is a $Z$-set in $M'$,
and $f(X \setminus A) \subset M'$. Take $\UUu \in \cov(M)$ with $B(f\bigr|_{X \setminus A},\UUu) \subset G$ and let
$\VVv' \in \cov(M')$ be a star refinement of the cover $\{U \cap M'\dd\ U \in \UUu\}$ of $M'$. Now by \THM{Z-maps},
there is a map $h'\dd X \to Y$ such that
\begin{equation}\label{eqn:aux7}
h'\bigr|_A = f\bigr|_A, \qquad h'(X \setminus A) \subset M' \setminus B'
\end{equation}
and $h'\bigr|_{X \setminus A}$ and $f\bigr|_{X \setminus A}$ are $\VVv'$-close. Fix metrics $d$ and $\varrho$
on $X$ and $Y$, respectively, such that for all $x,y \in X$:
\begin{equation}\label{eqn:nonexp}
\varrho(h'(x),h'(y)) \leqsl d(x,y).
\end{equation}
Further, we define $\lambda \in \CCc(X,[0,+\infty))$ and closed subsets $A_1,A_2,\ldots$ of $X$ as follows.
If $\partial A = \varempty$, put $\lambda \equiv 1$ and $A_n = X \setminus A\ (n \geqsl 1)$.
Otherwise, let $\lambda(x) = \dist_d(x,\partial A)$ and $A_n = \lambda^{-1}([\frac1n,+\infty)) \setminus A$.
Note that
\begin{equation}\label{eqn:incr}
A_k \subset \intt A_{k+1}\ (k \geqsl 1), \qquad \bigcup_{n=1}^{\infty} A_n = X \setminus A.
\end{equation}
Again by \THM{Z-maps}, there exists a map $h\dd X \to Y$ for which $h\bigr|_A = h'\bigr|_A$, $h(X \setminus A)
\subset M' \setminus (B' \cup B)$, $h\bigr|_{A_n}$ is a $Z$-embedding into $M'$ for all $n \in \NNN$,
$\varrho(h(x),h'(x)) < \lambda(x)$ for any $x \in X \setminus A$ and $h\bigr|_{X \setminus A}$ is $\VVv'$-close
to $h'\bigr|_{X \setminus A}$. We easily get (E1). We infer from the connections $h(X \setminus A) \subset
M' \setminus B'$ and $h\bigr|_A = h'\bigr|_A = f\bigr|_A$ that
\begin{equation}\label{eqn:aux9}
h(X \setminus A) \cap \overline{h(\partial A)} = \varempty.
\end{equation}
Note that if $\partial A = \varempty$, all the conditions (E1)--(E3) are basicly fulfilled. Therefore from now on,
we assume that the boundary of $A$ is nonempty.\par
Let $x_n,x \in X \setminus A$ be such that $h(x_n) \to h(x)\ (n\to\infty)$. We claim that there is $j \geqsl 1$ such
that $x_n \in A_j$ for almost all $n$. Suppose, for the contrary, that there is a subseqeunce $(y_n)_n$ of $(x_n)_n$
such that $y_n \notin A_n$. This says that $\lambda(y_n) < \frac1n$ and thus there is $a_n \in \partial A$ for which
\begin{equation}\label{eqn:aux11}
d(y_n,a_n) \to 0\ (n\to\infty).
\end{equation}
But then $\varrho(h(y_n),h'(y_n)) \leqsl \lambda(y_n) \to 0$, so $h'(y_n) \to h(x)$. What is more,
$\varrho(h'(y_n),h'(a_n)) \leqsl d(y_n,a_n) \to 0$ (by \eqref{eqn:nonexp}). This yields $h(a_n) = h'(a_n) \to h(x)$
and therefore $h(x) \in \overline{h(\partial A)}$, which denies \eqref{eqn:aux9} (since $x \in X \setminus A$).
So, there is $j \geqsl 1$ such that $x_n \in A_j$ for almost all $n$. But then $x_n \to x$, because of the facts that
$h\bigr|_{A_j}$ is a closed embedding of $A_j$ into $M'$, $h(x) \in M'$ and $h\bigr|_{X \setminus A}$ is one-to-one
(thanks to \eqref{eqn:incr}). We have shown that $h\bigr|_{X \setminus A}$ is an embedding. Similarly,
under the assumptions of (E3), $h\bigr|_{\overline{X \setminus A}}$ is an embedding. Indeed, let $y_n \in X
\setminus A$, $a \in \partial A$ and $h(y_n) \to h(a)$. By (E1), \eqref{eqn:aux9}, the assumptions of (E3) and thanks
to the closedness of $h(A_j)$ in $M'$, $y_n \notin A_j$ for almost all $n$. This implies that $\lambda(y_n) \to 0$
and there is a sequence $(a_n)_n$ of elements of $\partial A$ for which \eqref{eqn:aux11} is fulfilled.
As in the previous part of the proof, we show that $h'(y_n) \to h(a)$ and thus $h'(a_n) \to h(a) = h'(a)$.
This, combined with \eqref{eqn:aux7} and the assumptions of (E3), gives $a_n \to a$ and therefore $y_n \to a$,
by \eqref{eqn:aux11}.\par
It suffices to prove that $h(X \setminus A)$ is locally closed in $M$. Since $F := h(X \setminus A) \subset M'$
and $M'$ is open in $M$, it is enough to show that $F$ is locally closed in $M'$. If $y \in F$, then there is
$n \geqsl 1$ and $x \in \intt A_n$ for which $h(x) = y$. Since $h\bigr|_{X \setminus A}$ is an embedding, there
is an open in $M'$ set $V \subset M'$ such that $h(\intt A_n) = V \cap F$. But then $y \in V \cap F = V \cap h(A_n)$
and the latter set is closed in $V$.
\end{proof}

Under the notation of the statement of \PRO{emb}, the most interesting case of this result appears when
$K = \varempty$ or $K = \overline{f(\partial A)}$:

\begin{corr}{main}{cf. \mbox{\cite[Theorem~3.1]{and-mcch}}}
Let $X$ be completely metrizable of weight no greater than $\alpha$ and $A$ its closed subset. Let $M$ be
an $\alpha$-manifold and $B$ its $\sigma$-$Z$-set. Let $f \in \CCc(X,M)$ and $\UUu \in \cov(M)$.
If $f(X \setminus A) \cap \overline{f(\partial A)} = \varempty$ or $\overline{f(\partial A)}$ is a $Z$-set in $M$,
then there is a map $h\dd X \to M$ $\UUu$-homotopic to $f$ such that $h\bigr|_A = f\bigr|_A$,
$h\bigr|_{X \setminus A}$ is a $\sigma$-$Z$-embedding whose image is locally closed in $M$ and disjoint from $B$
and $\overline{h(X \setminus A)} = h(X \setminus A) \cup \overline{h(\partial A)}$. What is more,
if $f\bigr|_{\partial A}$ is an embedding and $\overline{f(\partial A)}$ \textup{[}$f(\partial A)$\textup{]}
is a $Z$-set in $M$, then $h$ is an embedding \textup{[}a $Z$-embedding\textup{]}.
\end{corr}

Observe that if the image of a map $f\dd X \to M$ is contained in a locally closed $\sigma$-$Z$-set $L$ of $M$, then
$K = \bar{L} \setminus L$ is closed, $f(X) \cap K = \varempty$ and $\overline{f(X)} \setminus K$ is a $\sigma$-$Z$-set
in $M$. This notice, \PRO{emb} and \COR{main} give (cf. \cite[Lemma~2.4]{and-mcch}):

\begin{cor}{homo}
Let $X$ be completely metrizable, $w(X) \leqsl \alpha$, $M$ an $\alpha$-manifold and $B$ its $\sigma$-$Z$-set.
Let $\UUu \in \cov(M)$.
\begin{enumerate}[\upshape(A)]
\item If $F\dd X \times I \to M$ is a homotopy such that the closure of $F(X \times \{0,1\})$ is a $Z$-set in $M$,
   then there is a homotopy $H\dd X \times I \to M$ $\UUu$-close to $F$ such that $H(\cdot,t) = F(\cdot,t)$
   for $t = 0,1$ and $H\bigr|_{X \times (0,1)}$ is a $\sigma$-$Z$-embedding whose image is disjoint from $B$.
   If, in addition, $F\bigr|_{X \times \{0,1\}}$ is an embedding, so is $H$.
\item If $f\dd X \to M$ is a map whose image is contained in a locally closed $\sigma$-$Z$-set in $M$, then there is
   a $\UUu$-homotopy $F\dd X \times I \to M$ such that $F(\cdot,0) = f$ and $F\bigr|_{X \times (0,1]}$ is
   a $\sigma$-$Z$-embedding whose image is locally closed and disjoint from $B$. If, in addition, $f$ is an embedding,
   so is $F$.
\end{enumerate}
\end{cor}

Our last goal is to give a sufficient condition under which a map can be extended by an open embedding.

\begin{thm}{open}
Let $A$ be a closed subset of $X$ such that $X \setminus A$ is an $\alpha$-manifold which has no more than $\alpha$
components (i.e. $w(X \setminus A) = \alpha$). Let $M$ be a subset of $Y$ such that $M$ is an $\alpha$-manifold
and let $B$ be a $Z$-set in $M$. Let $f\dd X \to Y$ be a map such that $f(X \setminus A) \subset M$. If $f$ is approximable
by maps $g\dd X \to Y$ such that $g\bigr|_A = f\bigr|_A$, $g(X \setminus A) \subset M$ and $g(X \setminus A) \cap
\overline{g(\partial A)} = \varempty$, then for each neighbourhood $G$ of $f\bigr|_{X \setminus A}$
in $\CCc(X \setminus A,M)$ there is a map $h\dd X \to Y$ satisfying the following conditions:
\begin{enumerate}[\upshape(O1)]
\item $h\bigr|_A = f\bigr|_A$ and $h(X \setminus A) \subset M \setminus B$,
\item $h\bigr|_{X \setminus A}$ is an open embedding,
\item $h\bigr|_{X \setminus A} \in G$.
\end{enumerate}
\end{thm}
\begin{proof}
Thanks to the assumption, we may and do assume that $f(X \setminus A) \cap \overline{f(\partial A)} = \varempty$.
Put $M' = M \setminus \overline{f(\partial A)}$ and $B' = B \cap M'$. Then $M'$ is an $\alpha$-manifold, $B'$ is a $Z$-set
in $M'$ and $f(X \setminus A) \subset M'$. Let $\SSs = \Emb^o(X \setminus A,M') \cap
\SsS_{M'}(X \setminus A,M' \setminus B')$. Since $\SsS_{M'}(X \setminus A,M' \setminus B')$ is open and dense
in $\CCc(X \setminus A,M')$, thus $\SSs$ is dense as well (by \THM{open-emb-dense}). Now the assertion follows
from \LEM{main2}.
\end{proof}

\begin{cor}{open}
Let $M$ and $N$ be two $\alpha$-manifolds such that $w(M) = \alpha$.
\begin{enumerate}[\upshape(A)]
\item If $F\dd M \times I \to N$ is a homotopy such that the closure of $F(X \times \{0,1\})$ is a $Z$-set in $N$,
   then for every $\UUu \in \cov(N)$ there is a $\UUu$-close to $F$ homotopy $H\dd M \times I \to N$ such that
   $H(\cdot,t) = F(\cdot,t)$ for $t=0,1$ and $H\bigr|_{M \times (0,1)}$ is an open embedding.
\item If $f\dd M \to N$ is such a map that $\overline{\im} f$ is a $Z$-set in $N$, then for every $\UUu \in \cov(N)$
   there is a $\UUu$-homotopy $F\dd M \times I \to N$ such that $F(\cdot,0) = f$ and $F\bigr|_{M \times (0,1]}$
   is an open embedding.
\end{enumerate}
\end{cor}

We end the paper with

\begin{exm}{false}
Let $f_t = \id_M\dd M \to M\ (t \in I)$, where $M$ is a manifold. By \LEM{disjoint}, there is no homotopy $(h_t)_{t \in I}$
such that $h_0 = h_1 = \id_M$ and $h\bigr|_{M \times (0,1)}$ is an embedding. This shows that the paper of West \cite{west}
contains an oversight. However, it was applied by Anderson and McCharen \cite{and-mcch}
once---in \cite[Lemma~2.4]{and-mcch}. Fortunately, the assertion of that lemma is true, which follows from our \COR{homo}.
\end{exm}

\end{document}